\documentclass{article}

\usepackage{times}
\usepackage{graphicx} 
\usepackage{subcaption}

\usepackage{natbib}

\usepackage{algorithm}
\usepackage{algorithmic}

\usepackage{amsmath}
\usepackage{amsfonts}
\usepackage[hyphens]{url}
\usepackage{graphicx}
\usepackage[utf8]{inputenc} 
\usepackage[T1]{fontenc}    
\usepackage{hyperref}       
\usepackage{url}            
\usepackage{booktabs}       
\usepackage{amsfonts}       
\usepackage{nicefrac}       
\usepackage{microtype}      
\usepackage{color}
\usepackage{amsthm}

\newcommand{\mbf}[1]{\mathbf{#1}}
\newcommand{\lhat}{\widehat{\mbf{L}}}
\newcommand{\uhat}{\widehat{\mbf{U}}}
\newtheorem{theorem}{Theorem}
\newtheorem{definition}{Definition}

\newtheorem{innercustomthm}{Theorem}
\newenvironment{customthm}[1]
  {\renewcommand\theinnercustomthm{#1}\innercustomthm}
  {\endinnercustomthm}

\usepackage{hyperref}

\usepackage[accepted]{icml2017}

\icmltitlerunning{An Efficient Algorithm for Low-Rank Approximation}

\begin{document} 

\twocolumn[
\icmltitle{An Efficient, Sparsity-Preserving, Online Algorithm for Low-Rank Approximation}



\icmlsetsymbol{equal}{*}

\begin{icmlauthorlist}
\icmlauthor{David Anderson}{equal,to}
\icmlauthor{Ming Gu}{equal,to}
\end{icmlauthorlist}

\icmlaffiliation{to}{University of California, Berkeley}

\icmlcorrespondingauthor{David Anderson}{davidanderson@berkeley.edu}
\icmlcorrespondingauthor{Ming Gu}{mgu@berkeley.edu}

\icmlkeywords{boring formatting information, machine learning, ICML}

\vskip 0.3in
]



\printAffiliationsAndNotice{\icmlEqualContribution} 

\begin{abstract} 
Low-rank  matrix  approximation  is  a fundamental tool in data analysis for processing large datasets, reducing noise, and finding important signals.
In this work, we present a novel truncated LU factorization called {\bf Spectrum-Revealing LU} (SRLU) for effective low-rank matrix approximation, and develop a fast algorithm to compute an SRLU factorization. We provide both matrix and singular value approximation error bounds for the SRLU approximation computed by our algorithm. Our analysis suggests that SRLU is competitive with the best low-rank matrix approximation methods, deterministic or randomized, in both computational complexity and approximation quality. 
Numeric experiments illustrate that SRLU preserves sparsity, highlights important data features and variables, can be efficiently updated, and calculates data approximations nearly as accurately as possible. 
To the best of our knowledge this is the first practical variant of the LU factorization for effective  and efficient low-rank matrix approximation. 
\end{abstract} 

\section{Introduction}

Low-rank approximation is an essential data processing technique for understanding  large or noisy data in diverse areas including data compression, image and pattern recognition, signal processing, compressed sensing, latent semantic indexing, anomaly detection, and recommendation systems.  Recent machine learning applications include training neural networks~\cite{journals/corr/JaderbergVZ14, journals/corr/KirkpatrickPRVD16}, second order online learning~\cite{conf/nips/LuoACL16}, representation learning~\cite{conf/naacl/WangMRS16}, and reinforcement learning~\cite{conf/nips/GhavamzadehLMM10}.  Additionally, a recent trend in machine learning is to include an approximation of second order information for better accuracy and faster convergence~\cite{conf/nips/KrummenacherMKB16}. 

In this work, we introduce a novel low-rank approximation algorithm called Spectrum-Revealing LU (SRLU) that can be efficiently computed and updated.  Furthermore, SRLU preserves sparsity and can identify important data variables and observations.  Our algorithm works on any data matrix, and achieves an approximation accuracy that only differs from the accuracy of the best approximation possible for any given rank by a constant factor.\footnote{The truncated SVD is known to provide the best low-rank matrix approximation, but it is rarely used for large scale practical data analysis. See a brief discussion of the SVD in supplemental material.}

The major innovation in SRLU is the efficient calculation of a truncated LU factorization of the form
\footnotesize
\begin{eqnarray}
\Pi_1\mbf{A}\Pi_2^T &=&\bordermatrix{ & \text{\tiny $k$} & \text{\tiny $m-k$} \cr \text{\tiny $k$} & \mbf{L}_{11}& \cr \text{\tiny $m-k$} & \mbf{L}_{21}&\mbf{I}_{n-k}}\bordermatrix{ & \text{\tiny $k$} & \text{\tiny $n-k$} \cr  & \mbf{U}_{11} & \mbf{U}_{12} \cr  &&\mbf{S}}\nonumber\\
&\approx& \begin{pmatrix}\mbf{L}_{11}\\ \mbf{L}_{21}\end{pmatrix}\begin{pmatrix}\mbf{U}_{11}& \mbf{U}_{12}\end{pmatrix}\nonumber\\
&\stackrel{\text{def}}{=}& \lhat\uhat,\nonumber
\end{eqnarray}
\normalsize
where $\Pi_1$ and $\Pi_2$ are judiciously chosen permutation matrices. The LU factorization is unstable, and in practice is implemented by pivoting (interchanging) rows during factorization, i.e. choosing permutation matrix $\Pi_1$.  For the truncated LU factorization to have any significance, nevertheless, complete pivoting (interchanging rows and columns) is necessary to guarantee that the factors $\widehat{\mbf{L}}$ and $\widehat{\mbf{U}}$ are well-defined and that their product accurately represents the original data.  Previously, complete pivoting was impractical as a matrix factorization technique because it requires accessing the entire data matrix at every iteration, but SRLU efficiently achieves complete pivoting through randomization and includes a deterministic follow-up procedure to ensure a hight quality low-rank matrix approximation, as supported by rigorous theory and numeric experiments.

\subsection{Background on the LU factorization}

Algorithm~\ref{alg:rightlooking_lu}  presents a basic implementation of the LU factorization, where the result is stored in place such that the upper triangular part of $\mbf{A}$ becomes $\mbf{U}$ and the strictly lower triangular part becomes the strictly lower part of $\mbf{L}$, with the diagonal of $\mbf{L}$ implicitly known to contain all ones. LU with partial pivoting finds the largest entry in the $i^{\text{th}}$ column from row $i$ to $m$ and pivots the row with that entry to the $i^{\text{th}}$ row.  LU with complete pivoting finds the largest entry in the submatrix $\mbf{A}_{i+1:m, i+1:n}$ and pivots that entry to $\mbf{A}_{i,i}$.  It is generally known and accepted that partial pivoting is sufficient for general, real-world data matrices in the context of linear equation solving.

\begin{algorithm}[H]
    \caption{The LU factorization}
    \label{alg:rightlooking_lu}
    \begin{algorithmic}[1]
    \STATE {\bf Inputs:} Data matrix $\mbf{A}\in\mathbb{R}^{m\times n}$
    \FOR{$i=1,2,\cdots,\min(m,n)$}
        \STATE Perform row and/or column pivots
        \FOR{$k=i+1,\cdots,m$}
            \STATE $\mbf{A}_{k,i} = \mbf{A}_{k,i} / \mbf{A}_{i,i}$
        \ENDFOR
        \STATE $\mbf{A}_{i+1:m, i+1:n} \; \mathrel{-}= \mbf{A}_{i+1:m,1:i} \cdot \mbf{A}_{1:i,i+1:n}$\label{alglineshur}
   \ENDFOR
   \end{algorithmic}
 \end{algorithm}
 
 \begin{algorithm}[H]
    \caption{Crout LU}
    \label{alg:crout_lu}
    \begin{algorithmic}[1]
   \STATE {\bfseries Inputs:} Data matrix $\mbf{A}\in\mathbb{R}^{m\times n}$, block size $b$
   \FOR{$j=0,b,2b,\cdots,\min(m,n)/b-1$}
        \STATE Perform column pivots
        \STATE $\mbf{A}_{j+1:m,j+1:j+b} -= $
        \STATE \quad $\mbf{A}_{j+1:m,1:j} \cdot \mbf{A}_{1:j,j+1:j+b}$. 
        \STATE Apply Algorithm~\ref{alg:rightlooking_lu} on $\mbf{A}_{j+1:m,j+1:j+b}$
        \STATE Apply the row pivots to other columns of $\mbf{A}$
        \STATE $\mbf{A}_{j+1:j+b,j+b+1:n} \; \mathrel{-}= $
        \STATE \quad $\mbf{A}_{j+1:j+b,1:j} \cdot \mbf{A}_{1:j,j+b+1:n}$
   \ENDFOR
  \end{algorithmic}
\end{algorithm}

Line 7 of Algorithm \ref{alg:rightlooking_lu} is known as the Schur update.  Given a sparse input, this is the only step of the LU factorization that causes fill.  As the algorithm progresses, fill will compound and may become dense, but the LU factorization, and truncated LU in particular, generally preserves some, if not most, of the sparsity of a sparse input.  A numeric illustration is presented below.

There are many variations of the LU factorization.  In Algorithm~\ref{alg:crout_lu} the Crout version of LU is presented in block form.  The column pivoting entails selecting the next $b$ columns so that the in-place LU step is performed on a non-singular matrix (provided the remaining entries are not all zero).  Note that the matrix multiplication steps are the bottleneck of this algorithm, requiring $O(mnb)$ operations each in general.

The LU factorization has been studied extensively since long before the invention of computers, with notable results from many mathematicians, including Gauss, Turing, and Wilkinson.  Current research on LU factorizations includes communication-avoiding implementations, such as tournament pivoting \cite{journals/siammax/KhabouDGG13}, sparse implementations \cite{journals/siamsc/GrigoriDL07}, and new computation of preconditioners \cite{journals/siamsc/ChowP15}.  A randomized approach to efficiently compute the LU factorization with complete pivoting recently appeared in \cite{journals/corr/MelgaardG15}.  These results are all in the context of linear equation solving, either directly or indirectly through an incomplete factorization used to precondition an iterative method.  This work repurposes the LU factorization to create a novel efficient and effective low-rank approximation algorithm using modern randomization technology.

\section{Previous Work}


\subsection{Low-Rank Matrix Approximation (LRMA)}
Previous work on low-rank data approximation includes the Interpolative Decomposition (ID) \cite{journals/siamsc/ChengGMR05}, the truncated QR with column pivoting factorization \cite{gu96}, and other deterministic column selection algorithms, such as in~\cite{batson12}.

Randomized algorithms have grown in popularity in recent years because of their ability to efficiently process large data matrices and because they can be supported with rigorous theory.  Randomized low-rank approximation algorithms generally fall into one of two categories: sampling algorithms and black box algorithms.  Sampling algorithms form data approximations from a random selection of rows and/or columns of the data.  Examples include \cite{journals/toc/DeshpandeRVW06, conf/approx/DeshpandeV06, journals/jacm/FriezeKV04, mahoney2009matrix}.  \cite{journals/siammax/DrineasMM08} showed that for a given approximate rank $k$, a randomly drawn subset $\mbf{C}$ of $c=O\left(k\log(k)\epsilon^{-2}\log\left(1/\delta\right)\right)$ columns of the data, a randomly drawn subset $\mbf{R}$ of $r=O\left(c\log(c)\epsilon^{-2}\log\left(1/\delta\right)\right)$ rows of the data, and setting $\mbf{U}=\mbf{C}^\dagger\mbf{A}\mbf{R}^\dagger$, then the matrix approximation error 
$\|\mbf{A}-\mbf{CUR}\|_F$ is at most a factor of $1+\epsilon$ from the optimal rank $k$ approximation with probability at least $1-\delta$.  Black box algorithms typically approximate a data matrix in the form
\begin{eqnarray}
\mbf{A} \approx \mbf{Q}^T\mbf{QA},\nonumber
\end{eqnarray}
where $\mbf{Q}$ is an orthonormal basis of the random projection (usually using SVD, QR, or ID).  The result of \cite{jl1984} provided the theoretical groundwork for these algorithms, which have been extensively studied \cite{journals/corr/abs-1207-6365, DBLP:journals/siamrev/HalkoMT11, martinsson06, journals/jcss/PapadimitriouRTV00, conf/focs/Sarlos06, woolfe2008fast, liberty2007randomized,journals/siamsc/Gu15}.  Note that the projection of an $m$-by-$n$ data matrix is of size $m$-by-$\ell$, for some oversampling parameter $\ell\ge k$, and $k$ is the target rank.  Thus the computational challenge is the orthogonalization of the projection (the random projection can be applied quickly, as described in these works).  A previous result on randomized LU factorizations for low-rank approximation was presented in \cite{journals/corr/AizenbudSA16}, but is uncompetitive in terms of theoretical results and computational performance with the work presented here.

For both sampling and black box algorithms the tuning parameter $\epsilon$ cannot be arbitrarily small, as the methods become meaningless if the number of rows and columns sampled (in the case of sampling algorithms) or the size of the random projection (in the case of black box algorithms) surpasses the size of the data.  A common practice is  $\epsilon\approx\frac{1}{2}$.

\subsection{Guaranteeing Quality}
Rank-revealing algorithms \cite{chan87} are LRMA algorithms that guarantee the approximation is of high quality by also capturing the rank of the data within a tolerance (see supplementary materials for definitions).  These methods, nevertheless, attempt to build an important submatrix of the data, and do not directly compute a low-rank approximation.  Furthermore, they do not attempt to capture all positive singular values of the data.  \cite{miranian03} introduced a new type of high-quality LRMA algorithms that can capture all singular values of a data matrix within a tolerance, but requires extra computation to bound approximations of the left and right null spaces of the data matrix.  Rank-revealing algorithms in general are designed around a definition that is not specifically appropriate for LRMA.

A key advancement of this work is a new definition of high quality low-rank approximation:
\begin{definition}\label{def:sr}
A rank-$k$ truncated LU factorization is \textbf{spectrum-revealing} if
\begin{eqnarray}
\left\Vert\mbf{A}-\lhat\uhat\right\Vert_2 &\le& q_1(k,m,n)\sigma_{k+1}\left(\mbf{A}\right)\nonumber
\end{eqnarray}
and
\begin{eqnarray}
\sigma_j\left(\lhat\uhat\right) &\ge& \frac{\sigma_j\left(\mbf{A}\right)}{q_2(k,m,n)}\nonumber
\end{eqnarray}
for $1\le j \le k$ and $q_1(k,m,n)$ and $q_2(k,m,n)$ are bounded by a low degree polynomial in $k$, $m$, and $n$.
\end{definition}
Definition \ref{def:sr} has precisely what we desire in an LRMA, and no additional requirements.  The constants, $q_1(k,m,n)$ and $q_2(k,m,n)$ are at least $1$ for any rank-$k$ approximation by~\cite{eckart1936approximation}.  This work shows theoretically and numerically that 
our algorithm, SRLU, is spectrum-revealing in that it always finds such $q_1$ and $q_2$, often with $q_1, q_2 = O(1)$ in practice.

\subsection{Low-Rank and Other Approximations in Machine Learning}
Low-rank and other approximation algorithms have appeared recently in a variety of machine learning applications.  In~\cite{conf/nips/KrummenacherMKB16}, randomized low-rank approximation is applied directly to the adaptive optimization algorithm \textsc{AdaGrad} to incorporate variable dependence during optimization to approximate the full matrix version of \textsc{AdaGrad} with a significantly reduced computational complexity.  In~\cite{journals/corr/KirkpatrickPRVD16}, a diagonal approximation of the posterior distribution of previous data is utilized to alleviate catastrophic forgetting.



\section{Main Contribution: Spectrum-Revealing LU (SRLU)}

Our algorithm for computing SRLU is composed of two subroutines: partially factoring the data matrix with randomized complete pivoting (TRLUCP) and performing swaps to improve the quality of the approximation (SRP). The first provides an efficient algorithm for computing a truncated LU factorization, whereas the second ensures the resulting approximation is provably reliable.

\subsection{Truncated Randomized LU with Complete Pivoting (TRLUCP)}\label{Sec:TRLUCP}

\begin{algorithm}[tb]
   \caption{TRLUCP}
   \label{alg:trlucp}
\begin{algorithmic}[1]
   \STATE {\bfseries Inputs:} Data matrix $\mbf{A}\in\mathbb{R}^{m\times n}$, target rank $k$, block size $b$, oversampling parameter $p\ge b$, random Gaussian matrix $\Omega\in\mathbb{R}^{p\times m}$, $\lhat$ and $\uhat$ are initially 0 matrices
   \STATE Calculate random projection $\mbf{R}=\Omega\mbf{A}\label{step:proj}$
   \FOR{$j=0,b,2b,\cdots,k-b$}
        \STATE Perform column selection algorithm on $\mbf{R}$ and swap columns of $\mbf{A}$\label{step:colsel}
        \STATE Update block column of $\lhat$\label{step:schur1}
        \STATE Perform block LU with partial row pivoting and swap rows of $\mbf{A}$\label{step:lu}
        \STATE Update block row of $\uhat$\label{step:schur2}
        \STATE Update $\mbf{R}$\label{step:r}
   \ENDFOR
\end{algorithmic}
\end{algorithm}

Intuitively, TRLUCP performs deterministic LU with partial row pivoting for some initial data with permuted columns.  TRLUCP uses a random projection of the Schur complement to cheaply find and move forward columns that are more likely to be representative of the data.  To accomplish this, Algorithm \ref{alg:trlucp} performs an iteration of block LU factorization in a careful order that resembles Crout LU reduction.  The ordering is reasoned as follows: LU with partial row pivoting cannot be performed until the needed columns are selected, and so column selection must first occur at each iteration.  Once a block column is selected, a partial Schur update must be performed on that block column before proceeding.  At this point, an iteration of block LU with partial row pivoting can be performed on the current block.  Once the row pivoting is performed, a partial Schur update of the block of pivoted rows of $\mbf{U}$ can be performed, which completes the factorization up to rank $j+b$.  Finally, the projection matrix $\mbf{R}$ can be cheaply updated to prepare for the next iteration.  Note that any column selection method may be used when picking column pivots from $\mbf{R}$, such as QR with column pivoting, LU with row pivoting, or even this algorithm can again be run on the subproblem of column selection of $\mbf{R}$.  The flop count of TRLUCP is dominated by the three matrix multiplication steps (lines 2, 5, and 7).  The total number of flops is
\begin{eqnarray}
F^{\text{TRLUCP}} &=& 2pmn+(m+n)k^2+O\left(k(m+n)\right).\nonumber
\end{eqnarray}
Note the transparent constants, and, because matrix multiplication is the bottleneck, this algorithm can be implemented efficiently in terms of both computation as well as memory usage.  Because the output of TRLUCP is only written once, the total number of memory writes is $(m+n-k)k$.  Minimizing the number of data writes by only writing data once significantly improves efficiency because writing data is typically one of the slowest computational operations.  Also worth consideration is the simplicity of the LU decomposition, which only involves three types of operations: matrix multiply, scaling, and pivoting.  By contrast, state-of-the-art calculation of both the full and truncated SVD requires a more complex process of bidiagonalization.  The projection $\mbf{R}$ can be updated efficiently to become a random projection of the Schur complement for the next iteration.  This calculation involves the current progress of the LU factorization and the random matrix $\Omega$, and is described in detail in the appendix.

\subsection{Spectrum-Revealing Pivoting (SRP)}

TRLUCP produces high-quality data approximations for almost all data matrices, despite the lack of theoretical guarantees, but can miss important rows or columns of the data.  Next, we develop an efficient variant of the existing rank-revealing LU algorithms \cite{gu96, miranian03} to rapidly detect and, if necessary, correct any possible matrix approximation failures of TRLUCP. 


Intuitively, the quality of the factorization can be tested by searching for the next choice of pivot in the Schur complement if the factorization continued and determining if the addition of that element would significantly improve the approximation quality.  If so, then the row and column with this element should be included in the approximation and another row and column should be excluded to maintain rank.  Because TRLUCP does not provide an updated Schur complement, the largest element in the Schur complement can be approximated by finding the column of $\mbf{R}$ with largest norm, performing a Schur update of that column, and then picking the largest element in that column.  Let $\alpha$ be this element, and, without loss of generality, assume it is the first entry of the Schur complement.  Denote:
\begin{equation}
\Pi_1\mbf{A}\Pi_2^T	= \begin{pmatrix}\mbf{L}_{11} &  & \\ \ell^T & 1 & \\ \mbf{L}_{31} & & \mbf{I}\end{pmatrix}\begin{pmatrix}\mbf{U}_{11} & u & \mbf{U}_{13} \\ & \alpha & s_{12}^T \\ & s_{21} & \mbf{S}_{22}\end{pmatrix}.\label{eqn:lu}
\end{equation}

Next, we must find the row and column that should be replaced if the row and column containing $\alpha$ are important.  Note that the smallest entry of $\mbf{L}_{11}\mbf{U}_{11}$ may still lie in an important row and column, and so the largest element of the inverse should be examined instead.  Thus we propose defining
\begin{eqnarray}
\overline{\mbf{A}}_{11}	\stackrel{\text{def}}{=} \begin{pmatrix}\mbf{L}_{11} & \\ \ell^T&1\end{pmatrix}\begin{pmatrix}\mbf{U}_{11}&u\\ &\alpha\end{pmatrix}\nonumber
\end{eqnarray}
and testing
\begin{equation}\label{eqn:test}
\|\overline{\mbf{A}}_{11}^{-1}\|_{\max} \le \frac{f}{|\alpha|}
\end{equation}
for a tolerance parameter $f>1$ that provides a control of accuracy
versus the number of swaps needed.  Should the test fail, the row and
column containing $\alpha$ are swapped with the row and column
containing the largest element in $\overline{\mbf{A}}_{11}^{-1}$.
Note that this element may occur in the last row or last column of
$\overline{\mbf{A}}_{11}^{-1}$, indicating only a column swap or row
swap respectively is needed.  When the swaps are performed, the
factorization must be updated to maintain truncated LU form. We have
developed a variant of the LU updating algorithm of \cite{updatelu} to
efficiently update the SRLU factorization. 

SRP can be implemented efficiently: each swap requires at most
$O\left(k(m+n)\right)$ operations, and
$\|\overline{\mbf{A}}_{11}^{-1}\|_{\max}$ can be quickly and reliably
estimated using \cite{inversecalc}.  An argument similar to that used
in \cite{miranian03} shows that each swap will increase $\left|{\bf
  det}\left( \overline{\mbf{A}}_{11}\right)\right|$ by a factor at
least $f$, hence will never repeat. At termination, SRP will ensure a
partial LU factorization of the form (\ref{eqn:lu}) that satisfies
condition (\ref{eqn:test}). We will discuss spectrum-revealing
properties of this factorization in Section 4.2.

It is possible to derive theoretical upper bounds on the worst number
of swaps necessary in SRP, but in practice, this number is zero for
most matrices, and does not exceed $3-5$ in the most pathological data
matrix of dimension at most 1000 we can contrive.


\begin{algorithm}[tb]
   \caption{Spectrum-Revealing Pivoting (SRP)}
   \label{alg:srlu}
\begin{algorithmic}[1]
   \STATE {\bfseries Input:} Truncated LU factorization $\mbf{A}\approx\lhat\uhat$, tolerance $f > 1$
   \WHILE{$\|\overline{\mbf{A}}_{11}^{-1}\|_{\max}>\frac{f}{|\alpha|}$}
   \STATE Set $\alpha$ to be the largest element in $\mbf{S}$ (or find an approximate $\alpha$ using $\mbf{R}$)
   \STATE Swap row and column containing $\alpha$ with row and column of largest element in $\overline{\mbf{A}}_{11}^{-1}$
   \STATE Update truncated LU factorization
   \ENDWHILE
\end{algorithmic}
\end{algorithm}

SRLU can be used effectively to approximate second order information in machine learning.  SRLU can be used as a modification to \textsc{AdaGrad} in a manner similar to the low-rank approximation method in~\cite{conf/nips/KrummenacherMKB16}.  Applying the initialization technique in this work, SRLU would likely provide an efficient and accurate adaptive stochastic optimization algorithm.  SRLU can also become a full-rank approximation (low-rank plus diagonal) by adding a diagonal approximation of the Schur complement.  Such an approximation could be appropriate for improving memory in artificial intelligence, such as in~\cite{journals/corr/KirkpatrickPRVD16}.  SRLU is also a freestanding compression algorithm.


\subsection{The CUR Decomposition with LU}\label{section:cur}

A natural extension of truncated LU factorizations is a CUR-type decomposition for increased accuracy \cite{mahoney2009matrix}:
\begin{eqnarray}
\Pi_1\mbf{A}\Pi_2^T &\approx& \lhat\left(\lhat^\dagger\mbf{A}\uhat^\dagger\right)\uhat\ \ \stackrel{\text{def}}{=}\ \ \lhat\mbf{M}\uhat.\nonumber
\end{eqnarray}
As with standard CUR, the factors $\lhat$ and $\uhat$ retain (much of) the sparsity of the original data, while $\mbf{M}$ is a small, $k$-by-$k$ matrix.  The CUR decomposition can improve the accuracy of an SRLU with minimal extra needed memory.  Extra computational time, nevertheless, is needed to calculate $\mbf{M}$.  A more efficient, approximate CUR decomposition can be obtained by replacing $\mbf{A}$ with a high quality approximation (such as an SRLU factorization of high rank) in the calculation of $\mbf{M}$.

\subsection{The Online SRLU Factorization}

Given a factored data matrix $\mbf{A}\in\mathbb{R}^{m\times n}$ and new observations $\mbf{B}\Pi_2^T=\bordermatrix{ & \text{\tiny $k$} & \text{\tiny $m-k$} \cr & \mbf{B}_1& \mbf{B}_2}\in\mathbb{R}^{s\times m}$, an augmented LU decomposition takes the form
\begin{eqnarray}
\begin{pmatrix}\Pi_1\mbf{A}\Pi_2^T\\ \mbf{B}\Pi_2^T\end{pmatrix} &=& \begin{pmatrix}\mbf{L}_{11} && \\ \mbf{L}_{21} & \mbf{I} & \\ \mbf{L}_{31} && \mbf{I}\end{pmatrix}\begin{pmatrix}\mbf{U}_{11}&\mbf{U}_{12} \\ & \mbf{S} \\ & \mbf{S}^{\text{new}}\end{pmatrix},\nonumber
\end{eqnarray}
where $\mbf{L}_{31}=\mbf{B}_1\mbf{U}_{11}^{-1}$ and $\mbf{S}^{\text{new}}=\mbf{B}_2-\mbf{B}_1\mbf{U}_{11}^{-1}\mbf{U}_{12}$.  An SRLU factorization can then be obtained by simply performing correcting swaps.  For a rank-1 update, at most 1 swap is expected (although examples can be constructed that require more than one swap), which requires at most $O\left(k\left(m+n\right)\right)$ flops.  By contrast, the URV decomposition of \cite{journals/tsp/Stewart92} is $O\left(n^2\right)$, while SVD updating requires $O\left(\left(m+n\right)\min^2\left(m,n\right)\right)$ operations in general, or $O\left(\left(m+n\right)\min\left(m,n\right)\log_2^2\epsilon\right)$ for a numerical approximation with the fast multipole method.  

\section{Theoretical Results for SRLU Factorizations}\label{section:theory}

\subsection{Analysis of General Truncated LU Decompositions}

\begin{theorem}\label{theorem1}
Let $\left(\cdot\right)_s$ denote the rank-$s$ truncated SVD for $s\le k\ll m,n$.  Then for any truncated LU factorization with Schur complement $\mbf{S}$:
\begin{eqnarray}
\|\Pi_1\mbf{A}\Pi_2^T-\lhat\uhat\| &=& \|\mbf{S}\|\nonumber
\end{eqnarray}
for any norm, and
\begin{eqnarray}
\|\Pi_1\mbf{A}\Pi_2^T-\left(\lhat\uhat\right)_s\|_2  &\le& 2\|\mbf{S}\|_2+\sigma_{s+1}\left(\mbf{A}\right).\nonumber
\end{eqnarray}
\end{theorem}
\begin{theorem}\label{theoremb}
For a general rank-$k$ truncated LU decomposition, we have for all $1 \leq j \leq k$, 
\begin{eqnarray}
\sigma_j\left(\mbf{A}\right) \le \sigma_j\left(\lhat\uhat\right)\left(1+\left(1+\frac{\|\mbf{S}\|_2}{\sigma_k\left(\lhat\uhat\right)}\right)\frac{\|\mbf{S}\|_2}{\sigma_j\left(\mbf{A}\right)}\right).\nonumber
\end{eqnarray}
\end{theorem}
\begin{theorem}{CUR Error Bounds.}\label{theorem4}
\begin{eqnarray}
\|\Pi_1\mbf{A}\Pi_2^T-\lhat\mbf{M}\uhat\|_2 &\le& 2\|\mbf{S}\|_2\nonumber
\end{eqnarray}
and
\begin{eqnarray}
\|\Pi_1\mbf{A}\Pi_2^T-\lhat\mbf{M}\uhat\|_F &\le& \|\mbf{S}\|_F.\nonumber
\end{eqnarray}
\end{theorem}
Theorem \ref{theorem1} simply concludes that the approximation is accurate if the Schur complement is small, but the singular value bounds of Theorem \ref{theoremb} are needed to guarantee that the approximation retains structural properties of the original data, such as an accurate approximation of the rank and the spectrum.  Furthermore, singular values bounds can be significantly stronger than the more familiar norm error bounds that appear in Theorem \ref{theorem1}.  Theorem \ref{theoremb} provides a general framework for singular value bounds, and bounding the terms in this theorem provided guidance in the design and development of SRLU.  Just as in the case of deterministic LU with complete pivoting, the sizes of  $\frac{\|\mbf{S}\|_2}{\sigma_k\left(\lhat\uhat\right)}$ and $\frac{\|\mbf{S}\|_2}{\sigma_j\left(\lhat\uhat\right)}$ range from moderate to small for almost all data matrices of practical interest. They, nevertheless, cannot be effectively bounded for a general TRLUCP factorization, implying the need for Algorithm \ref{alg:srlu} to ensure that these terms are controlled.  While the error bounds in Theorem \ref{theorem4} for the CUR decomposition do not improve upon the result in Theorem \ref{theorem1}, CUR bounds for SRLU specifically will be considerably stronger.  Next, we present our main theoretical contributions.

\subsection{Analysis of the Spectrum-Revealing LU Decomposition}\label{section:srlu}
\begin{theorem}{(SRLU Error Bounds.)}\label{theoremq}
For $j\le k$ and $\gamma=O\left(fk\sqrt{mn}\right)$, SRP produces a rank-$k$ SRLU factorization with 
\begin{eqnarray}
&\|\Pi_1\mbf{A}\Pi_2^T-\lhat\uhat\|_2 \le \gamma\sigma_{k+1}\left(\mbf{A}\right),&\nonumber\\
&\|\Pi_1\mbf{A}\Pi_2^T-\left(\lhat\uhat\right)_j\|_2 \le \sigma_{j+1}\left(\mbf{A}\right)\left(1+2\gamma\frac{\sigma_{k+1}\left(\mbf{A}\right)}{\sigma_{j+1}\left(\mbf{A}\right)}\right)&\nonumber
\end{eqnarray}
\end{theorem}
Theorem \ref{theoremq} is a special case of Theorem \ref{theorem1} for SRLU factorizations.  For a data matrix with a rapidly decaying spectrum, the right-hand side of the second inequality is close to $\sigma_{j+1}\left(\mbf{A}\right)$, a substantial improvement over the sharpness of the bounds in \cite{journals/siammax/DrineasMM08}.

\begin{theorem}{(SRLU Spectral Bound)}.\label{thm:srlubnd}
For $1\le j\le k$, SRP produces a rank-$k$ SRLU factorization with 
\begin{eqnarray}
\frac{\sigma_j\left(\mbf{A}\right)}{1+\tau\frac{\sigma_{k+1}\left(\mbf{A}\right)}{\sigma_j\left(\mbf{A}\right)}} \le \sigma_j\left(\lhat\uhat\right) \le \sigma_j\left(\mbf{A}\right)\left(1+\tau\frac{\sigma_{k+1}\left(\mbf{A}\right)}{\sigma_j\left(\mbf{A}\right)}\right)\nonumber
\end{eqnarray}
for $\tau\le O\left(mnk^2f^3\right)$.
\end{theorem}
While the worst case upper bound on $\tau$ is large, it is dimension-dependent, and $j$ and $k$ may be chosen so that $\frac{\sigma_{k+1}\left(\mbf{A}\right)}{\sigma_j\left(\mbf{A}\right)}$ is arbitrarily small compared to $\tau$.  In particular, if $k$ is the numeric rank of $\mbf{A}$, then the singular values of the approximation are numerically equal to those of the data.

These bounds are {\bf problem-specific bounds} because their quality depends on the spectrum of the original data, rather than universal constants that appear in previous results.  The benefit of these problem-specific bounds is that an approximation of data with a rapidly decaying spectrum is guaranteed to be high-quality.  Furthermore, if $\sigma_{k+1}\left(\mbf{A}\right)$ is not small compared to $\sigma_j\left(\mbf{A}\right)$, then no high-quality low-rank approximation is possible in the 2 and Frobenius norms.  Thus, in this sense, the bounds presented in Theorems \ref{theoremq} and \ref{thm:srlubnd} are optimal.

Given a high-quality rank-$k$ truncated LU factorization, Theorem \ref{thm:srlubnd} ensures that a low-rank approximation of rank $\ell$ with $\ell<k$ of the compressed data is an accurate rank-$\ell$ approximation of the full data.  The proof of this theorem centers on bounding the terms in Theorems \ref{theorem1} and \ref{theoremb}.  Experiments will show that $\tau$ is small in almost all cases.

Stronger results are achieved with the CUR version of SRLU:
\begin{theorem}\label{theorem5}
\begin{eqnarray}
\|\Pi_1\mbf{A}\Pi_2^T-\lhat\mbf{M}\uhat\|_2 &\le& 2\gamma\sigma_{k+1}\left(\mbf{A}\right)\nonumber
\end{eqnarray}
and
\begin{eqnarray}
\|\Pi_1\mbf{A}\Pi_2^T-\lhat\mbf{M}\uhat\|_F &\le& \omega\sigma_{k+1}\left(\mbf{A}\right),\nonumber
\end{eqnarray}
where $\gamma=O\left(fk\sqrt{mn}\right)$ is the same as in Theorem \ref{theoremq}, and $\omega=O\left(fkmn\right)$.
\end{theorem}
\begin{theorem}\label{theorem:cur}
If $\sigma_j^2\left(\mbf{A}\right)>2\|\mbf{S}\|_2^2$ then
\begin{eqnarray}
\sigma_j\left(\mbf{A}\right) \ge\sigma_j\left(\lhat\mbf{M}\uhat\right) \ge \sigma_j\left(\mbf{A}\right)\sqrt{1-2\gamma\left(\frac{\sigma_{k+1}\left(\mbf{A}\right)}{\sigma_j\left(\mbf{A}\right)}\right)^2}\nonumber
\end{eqnarray}
for $\gamma=O\left(mnk^2f^2\right)$ and $f$ is an input parameter controlling a tradeoff of quality vs. speed as before.
\end{theorem}
As before, the constants are small in practice.  Observe that for most real data matrices, their singular values decay with increasing $j$. For such matrices this result is significantly stronger than Theorem \ref{thm:srlubnd}.

\section{Experiments}\label{section:experiments}

\subsection{Speed and Accuracy Tests}

\begin{figure*}[t]
\centering
	\begin{subfigure}[t]{0.4\textwidth}
		\centering
		\includegraphics[width=2.5in]{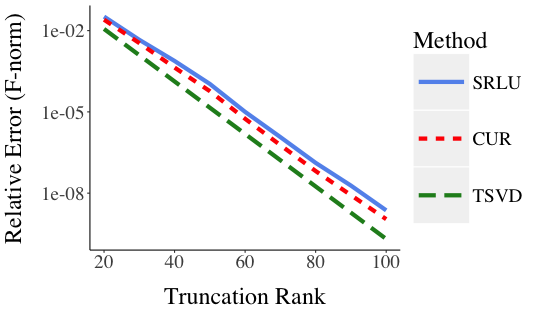}
		\caption{Spectral Decay = 0.8}
	\end{subfigure}%
	\begin{subfigure}[t]{0.4\textwidth}
		\centering
		\includegraphics[width=2.5in]{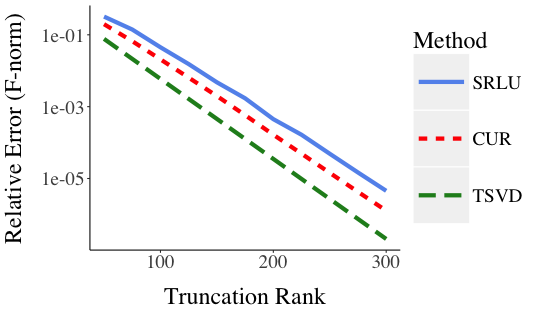}
		\caption{Spectral Decay = 0.95}
	\end{subfigure}
	\caption{Accuracy experiment on random 1000x1000 matrices with different rates of spectral decay.}
	\label{accuracy}
\end{figure*} 

\begin{figure*}[t]
\centering
	\begin{subfigure}[t]{0.4\textwidth}
		\centering
		\includegraphics[width=2.5in]{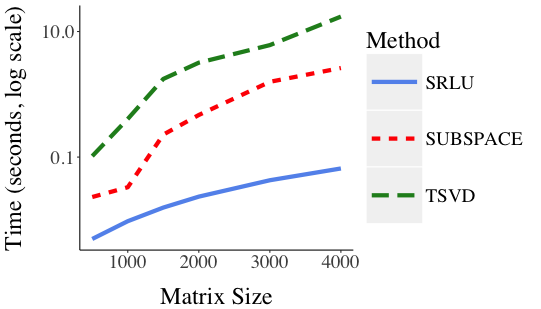}
		\caption{Rank-100 Factorizations}
	\end{subfigure}%
	\begin{subfigure}[t]{0.4\textwidth}
		\centering
		\includegraphics[width=2.5in]{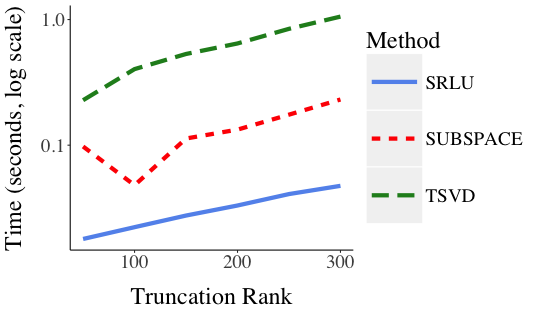}
		\caption{Time vs. Truncation Rank}
	\end{subfigure}
	\caption{Time experiment on various random matrices, and a time experiment on a 1000x1000 matrix with varying truncation ranks.}
	\label{time}
\end{figure*} 

In Figure \ref{accuracy}, the accuracy of our method is compared to the accuracy of the truncated SVD.  Note that SRLU did not perform any swaps in these experiments.  ``CUR" is the CUR version of the output of SRLU.  Note that both methods exhibits a convergence rate similar to that of the truncated SVD (TSVD), and so only a constant amount of extra work is needed to achieve the same accuracy.  When the singular values decay slowly, the CUR decomposition provides a greater accuracy boost.  In Figure \ref{time}, the runtime of SRLU is compared to that of the truncated SVD, as well as Subspace Iteration \cite{journals/siamsc/Gu15}.  Note that for Subspace Iteration, we choose iteration parameter $q=0$ and do not measure the time of applying the random projection, in acknowledgement that fast methods exist to apply a random projection to a data matrix.  Also, the block size implemented in SRLU is significantly smaller than the block size used by the standard software LAPACK, as the size of the block size affects the size of the projection.  See supplement for additional details.  All numeric experiments were run on NERSC's Edison.  For timing experiments, the truncated SVD is calculated with PROPACK.


Even more impressive, the factorization stage of SRLU becomes arbitrarily faster than the standard implementation of the LU decomposition.  Although the standard LU decomposition is not a low-rank approximation algorithm, it is known to be roughly 10 times faster than the SVD \cite{demmel97}.  See appendix for details.

Next, we compare SRLU against competing algorithms.  In \cite{conf/icml/UbaruMS15}, error-correcting codes are introduced to yield improved accuracy over existing random projection low-rank approximation algorithms.  Their algorithm, denoted {\tt Dual BCH}, is compared against SRLU as well as two other random projection methods: {\tt Gaus.}, which uses a Gaussian random projection, and {\tt SRFT}, which uses a Fourier transform to apply a random projection.  We test the spectral norm error of these algorithms on matrices from the sparse matrix collection in \cite{matrixcollection}.

\begin{table}[t]
\caption{Errors of low-rank approximations of the given target rank.  SRLU is measured using the CUR version of the factorization.}
\label{lucp-table}
\vskip 0.15in
\begin{center}
\begin{small}
\begin{sc}
\begin{tabular}{lccccc}
\hline
\abovespace\belowspace
Data & $k$ & {\tt Gaus.} & {\tt SRFT} & {\tt Dual BCH} & {\tt SRLU} \\
\hline
\abovespace
${\tt S80PI_n1}$  	&  63 	& 3.85 		& 3.80 		& 3.81 		& 2.84 \\
${\tt deter3}$	& 127 	& 9.27 		& 9.30 		& 9.26 		& 8.30\\
${\tt lc3d}$ 	& 63 		& 18.39 	& 16.36 	& 15.49 	& 16.94 \\
${\tt lc3d}$ 	& 78 		&  			&  			&  			& 15.11 \\
\hline
\end{tabular}
\end{sc}
\end{small}
\end{center}
\vskip -0.1in
\end{table}

In Table \ref{lucp-table}, results for SRLU are averaged over 5 experiments.  Using tuning parameter $f=5$, no swaps were needed in all cases.  The matrices being tested are sparse matrices from various engineering problems.  ${\tt S80PI_n1}$ is 4,028 by 4,028, {\tt deter3} is 7,647 by 21,777, and {\tt lp\_ceria3d} (abbreviated {\tt lc3d}) is 3,576 by 4,400.  Note that SRLU, a more efficient algorithm, provides a better approximation in two of the three experiments.  With a little extra oversampling, a practical assumption due to the speed advantage, SRLU achieves a competitive quality approximation.  The oversampling highlights an additional and unique advantage of SRLU over competing algorithms: if more accuracy is desired, then the factorization can simply continue as needed.

\subsection{Sparsity Preservation Tests}

\begin{figure}[h]
	\centering
	\begin{subfigure}[t]{0.5\textwidth}
		\centering
		\includegraphics[width=3.2in]{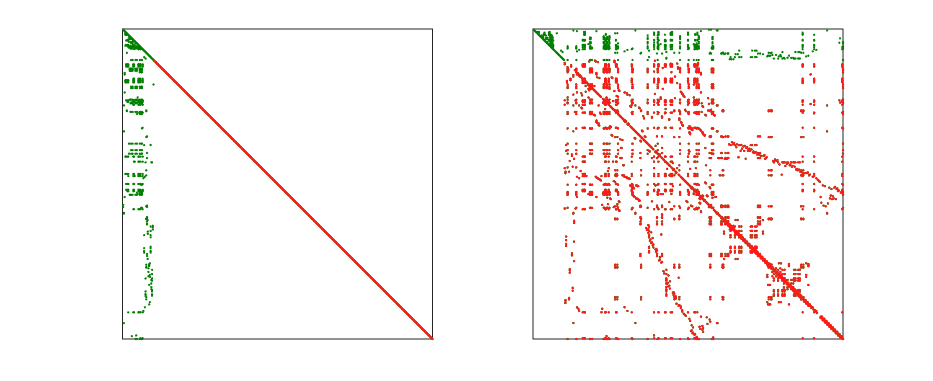}
		\caption{$\mbf{L}$ and $\mbf{U}$ patterns of a low-rank factorization}
	\end{subfigure}
	\begin{subfigure}[t]{0.5\textwidth}
		\centering
		\includegraphics[width=3.2in]{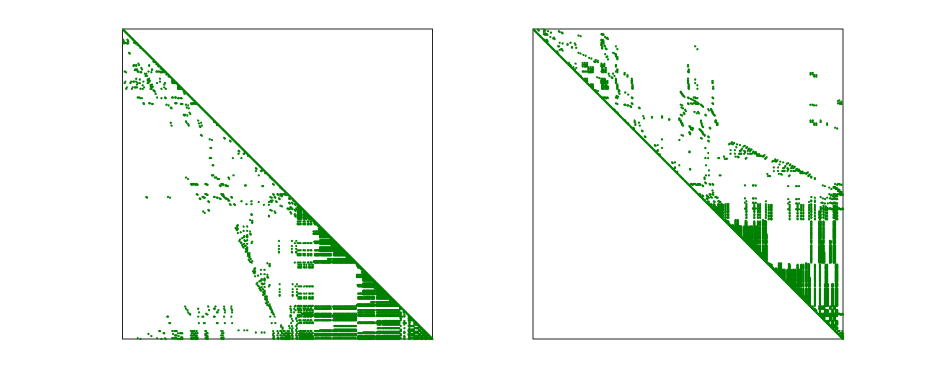}
		\caption{$\mbf{L}$ and $\mbf{U}$ patterns of the full factorization}
	\end{subfigure}
	\caption{The sparsity patterns of the L and U matrices of a rank 43 SRLU factorization, followed by the sparsity pattern of the L and U matrices of a full LU decomposition of the same data.  For the SRLU factorization, the green entries compose the low-rank approximation of the data.  The red entries are the additional data needed for an exact factorization.}
	\label{sparse_srlu}
\end{figure}

The SRLU factorization is tested on sparse, unsymmetric matrices from \cite{matrixcollection}.  Figure \ref{sparse_srlu} shows the sparsity patterns of the factors of an SRLU factorization of a sparse data matrix representing a circuit simulation (\texttt{oscil\_dcop}), as well as a full LU decomposition of the data.  Note that the LU decomposition preserves the sparsity of the data initially, but the full LU decomposition becomes dense.  Several more experiments are shown in the supplement.

\subsection{Towards Feature Selection}

\begin{figure}[ht!]
	\centering
	\begin{subfigure}[t]{0.15\textwidth}
		\centering
		\includegraphics[width=0.9in]{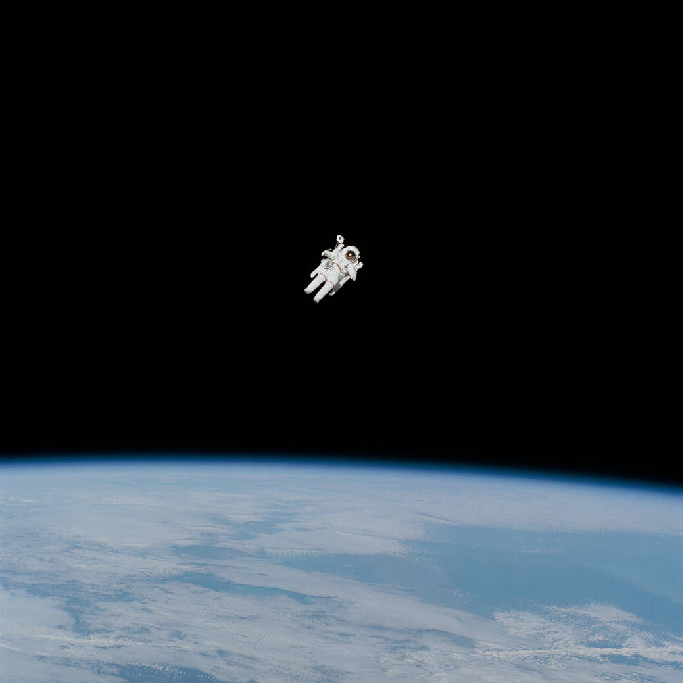}
		\caption{Original}
	\end{subfigure}%
	\begin{subfigure}[t]{0.15\textwidth}
		\centering
		\includegraphics[width=0.9in]{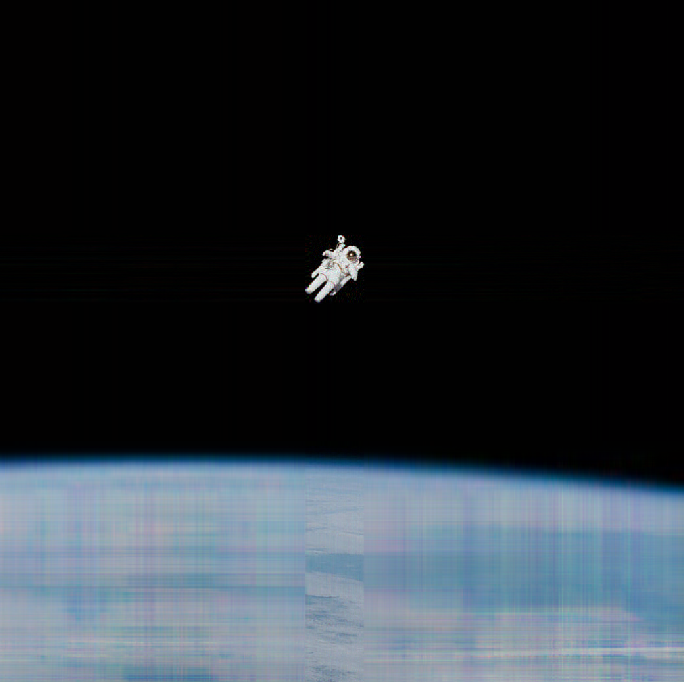}
		\caption{SRLU}
	\end{subfigure}%
	\begin{subfigure}[t]{0.15\textwidth}
		\centering
		\includegraphics[width=0.9in]{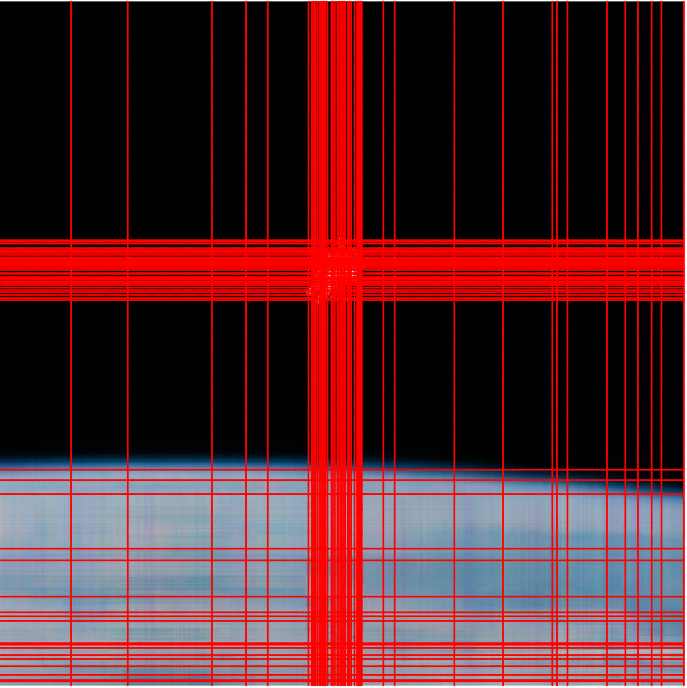}
		\caption{Rows and Cols.}
	\end{subfigure}
	\caption{Image processing example.  The original image \cite{nasa2jpg}, a rank-50 approximation with SRLU, and a highlight of the rows and columns selected by SRLU.}
	\label{astro_ex}
\end{figure}

An image processing example is now presented to illustrate the benefit of highlighting important rows and columns selection.  In Figure \ref{astro_ex} an image is compressed to a rank-50 approximation using SRLU.  Note that the rows and columns chosen overlap with the astronaut and the planet, implying that minimal storage is needed to capture the black background, which composes approximately two thirds of the image.  While this result cannot be called feature selection per se, the rows and columns selected highlight where to look for features: rows and/or columns are selected in a higher density around the astronaut, the curvature of the planet, and the storm front on the planet.

\subsection{Online Data Processing}

Online SRLU is tested here on the Enron email corpus \cite{Lichman:2013}.  The documents were initially reverse-sorted by the usage of the most common word, and then reverse-sorted by the second most, and this process was repeated for the five most common words (the top five words were used significantly more than any other), so that the most common words occurred most at the end of the corpus.  The data contains 39,861 documents and 28,102 words/terms, and an initial SRLU factorization of rank 20 was performed on the first 30K documents.  The initial factorization contained none of the top five words, but, after adding the remaining documents and updating, the top three were included in the approximation.  The fourth and fifth words `market' and `california' have high covariance with at least two of the three top words, and so their inclusion may be redundant in a low-rank approximation.

\section{Conclusion}
We have presented SRLU, a low-rank approximation method with many desirable properties: efficiency, accuracy, sparsity-preservation, the ability to be updated, and the ability to highlight important data features and variables.  Extensive theory and numeric experiments have illustrated the efficiency and effectiveness of this method.

\section*{Acknowledgements} 
This research was supported in part by NSF Award CCF-1319312.


\bibliography{lucp_icml}

\begin{thebibliography}{42}
\providecommand{\natexlab}[1]{#1}
\providecommand{\url}[1]{\texttt{#1}}
\expandafter\ifx\csname urlstyle\endcsname\relax
  \providecommand{\doi}[1]{doi: #1}\else
  \providecommand{\doi}{doi: \begingroup \urlstyle{rm}\Url}\fi

\bibitem[Aizenbud et~al.(2016)Aizenbud, Shabat, and
  Averbuch]{journals/corr/AizenbudSA16}
Aizenbud, Y., Shabat, G., and Averbuch, A.
\newblock Randomized lu decomposition using sparse projections.
\newblock \emph{CoRR}, abs/1601.04280, 2016.

\bibitem[Batson et~al.(2012)Batson, Spielman, and Srivastava]{batson12}
Batson, J., Spielman, D.~A., and Srivastava, N.
\newblock Twice-ramanujan sparsifiers.
\newblock \emph{SIAM Journal on Computing}, 41\penalty0 (6):\penalty0
  1704--1721, 2012.

\bibitem[Chan(1987)]{chan87}
Chan, T.~F.
\newblock Rank revealing qr factorizations.
\newblock \emph{Linear algebra and its applications}, \penalty0
  (88/89):\penalty0 67--82, 1987.

\bibitem[Cheng et~al.(2005)Cheng, Gimbutas, Martinsson, and
  Rokhlin]{journals/siamsc/ChengGMR05}
Cheng, H., Gimbutas, Z., Martinsson, P.-G., and Rokhlin, V.
\newblock On the compression of low rank matrices.
\newblock \emph{SIAM J. Scientific Computing}, 26\penalty0 (4):\penalty0
  1389--1404, 2005.

\bibitem[Chow \& Patel(2015)Chow and Patel]{journals/siamsc/ChowP15}
Chow, E. and Patel, A.
\newblock Fine-grained parallel incomplete lu factorization.
\newblock \emph{SIAM J. Scientific Computing}, 37\penalty0 (2), 2015.

\bibitem[Clarkson \& Woodruff(2012)Clarkson and
  Woodruff]{journals/corr/abs-1207-6365}
Clarkson, K.~L. and Woodruff, D.~P.
\newblock Low rank approximation and regression in input sparsity time.
\newblock \emph{CoRR}, abs/1207.6365, 2012.

\bibitem[David \& Hu(2011)David and Hu]{timdavis}
David, T.~A. and Hu, Y.
\newblock The university of florida sparse matrix collection.
\newblock \emph{ACM Transactions on Mathematical Software}, 38:\penalty0 1--25,
  2011.
\newblock URL \url{http://www.cise.ufl.edu/research/sparse/matrices}.

\bibitem[Davis \& Hu(2011)Davis and Hu]{matrixcollection}
Davis, T.~A. and Hu, Y.
\newblock The university of florida sparse matrix collection.
\newblock \emph{ACM Transactions on Mathematical Software}, 38:\penalty0
  1:1--1:25, 2011.
\newblock URL \url{http://www.cise.ufl.edu/research/sparse/matrices}.

\bibitem[Demmel(1997)]{demmel97}
Demmel, J.
\newblock \emph{Applied Numerical Linear Algebra}.
\newblock SIAM, 1997.

\bibitem[Deshpande \& Vempala(2006)Deshpande and
  Vempala]{conf/approx/DeshpandeV06}
Deshpande, A. and Vempala, S.
\newblock Adaptive sampling and fast low-rank matrix approximation.
\newblock In \emph{APPROX-RANDOM}, volume 4110 of \emph{Lecture Notes in
  Computer Science}, pp.\  292--303. Springer, 2006.

\bibitem[Deshpande et~al.(2006)Deshpande, Rademacher, Vempala, and
  Wang]{journals/toc/DeshpandeRVW06}
Deshpande, A., Rademacher, L., Vempala, S., and Wang, G.
\newblock Matrix approximation and projective clustering via volume sampling.
\newblock \emph{Theory of Computing}, 2\penalty0 (12):\penalty0 225--247, 2006.

\bibitem[Drineas et~al.(2008)Drineas, Mahoney, and
  Muthukrishnan]{journals/siammax/DrineasMM08}
Drineas, P., Mahoney, M.~W., and Muthukrishnan, S.
\newblock Relative-error cur matrix decompositions.
\newblock \emph{SIAM J. Matrix Analysis Applications}, 30\penalty0
  (2):\penalty0 844--881, 2008.

\bibitem[Eckart \& Young(1936)Eckart and Young]{eckart1936approximation}
Eckart, C. and Young, G.
\newblock The approximation of one matrix by another of lower rank.
\newblock \emph{Psychometrika}, 1\penalty0 (3):\penalty0 211--218, 1936.

\bibitem[Frieze et~al.(2004)Frieze, Kannan, and
  Vempala]{journals/jacm/FriezeKV04}
Frieze, A.~M., Kannan, R., and Vempala, S.
\newblock Fast monte-carlo algorithms for finding low-rank approximations.
\newblock \emph{J. ACM}, 51\penalty0 (6):\penalty0 1025--1041, 2004.

\bibitem[Ghavamzadeh et~al.(2010)Ghavamzadeh, Lazaric, Maillard, and
  Munos]{conf/nips/GhavamzadehLMM10}
Ghavamzadeh, M., Lazaric, A., Maillard, O.-A., and Munos, R.
\newblock Lstd with random projections.
\newblock In \emph{NIPS}, pp.\  721--729, 2010.

\bibitem[Golub \& van Loan(2013)Golub and van Loan]{golub13}
Golub, G.~H. and van Loan, C.~F.
\newblock \emph{Matrix Computations}.
\newblock JHU Press, 4th edition, 2013.

\bibitem[Gondzio(2007)]{updatelu}
Gondzio, J.
\newblock Stable algorithm for updating dense lu factorization after row or
  column exchange and row and column addition or deletion.
\newblock \emph{Optimization: A journal of Mathematical Programming and
  Operations Research}, pp.\  7--26, 2007.

\bibitem[Grigori et~al.(2007)Grigori, Demmel, and
  Li]{journals/siamsc/GrigoriDL07}
Grigori, L., Demmel, J., and Li, X.~S.
\newblock Parallel symbolic factorization for sparse lu with static pivoting.
\newblock \emph{SIAM J. Scientific Computing}, 29\penalty0 (3):\penalty0
  1289--1314, 2007.

\bibitem[Gu(2015)]{journals/siamsc/Gu15}
Gu, M.
\newblock Subspace iteration randomization and singular value problems.
\newblock \emph{SIAM J. Scientific Computing}, 37\penalty0 (3), 2015.

\bibitem[Gu \& Eisenstat(1996)Gu and Eisenstat]{gu96}
Gu, M. and Eisenstat, S.~C.
\newblock Efficient algorithms for computing a strong rank-revealing qr
  factorization.
\newblock \emph{SIAM J. Sci. Comput.}, 17\penalty0 (4):\penalty0 848--869,
  1996.

\bibitem[Halko et~al.(2011)Halko, Martinsson, and
  Tropp]{DBLP:journals/siamrev/HalkoMT11}
Halko, N., Martinsson, P.-G., and Tropp, J.~A.
\newblock Finding structure with randomness: Probabilistic algorithms for
  constructing approximate matrix decompositions.
\newblock \emph{SIAM Review}, 53\penalty0 (2):\penalty0 217--288, 2011.

\bibitem[Higham \& Relton(2015)Higham and Relton]{inversecalc}
Higham, N.~J. and Relton, S.~D.
\newblock Estimating the largest entries of a matrix.
\newblock 2015.
\newblock URL \url{http://eprints.ma.man.ac.uk/}.

\bibitem[Jaderberg et~al.(2014)Jaderberg, Vedaldi, and
  Zisserman]{journals/corr/JaderbergVZ14}
Jaderberg, M., Vedaldi, A., and Zisserman, A.
\newblock Speeding up convolutional neural networks with low rank expansions.
\newblock \emph{CoRR}, abs/1405.3866, 2014.

\bibitem[Johnson \& Lindenstrauss(1984)Johnson and Lindenstrauss]{jl1984}
Johnson, W.~B. and Lindenstrauss, J.
\newblock Extensions of lipschitz mappings into a hilbert space.
\newblock \emph{Contemporary Mathematics}, 26:\penalty0 189--206, 1984.

\bibitem[Khabou et~al.(2013)Khabou, Demmel, Grigori, and
  Gu]{journals/siammax/KhabouDGG13}
Khabou, A., Demmel, J., Grigori, L., and Gu, M.
\newblock Lu factorization with panel rank revealing pivoting and its
  communication avoiding version.
\newblock \emph{SIAM J. Matrix Analysis Applications}, 34\penalty0
  (3):\penalty0 1401--1429, 2013.

\bibitem[Kirkpatrick et~al.(2017)Kirkpatrick, Pascanu, Rabinowitz, Veness,
  Desjardins, Rusu, Milan, Quan, Ramalho, Grabska-Barwinska, Hassabis, Clopath,
  Kumaran, and Hadsell]{journals/corr/KirkpatrickPRVD16}
Kirkpatrick, J., Pascanu, R., Rabinowitz, N.~C., Veness, J., Desjardins, G.,
  Rusu, A.~A., Milan, K., Quan, J., Ramalho, T., Grabska-Barwinska, A.,
  Hassabis, D., Clopath, C., Kumaran, D., and Hadsell, R.
\newblock Overcoming catastrophic forgetting in neural networks.
\newblock \emph{Proceedings of the National Academy of Sciences}, 114\penalty0
  (13):\penalty0 3521--3526, 2017.

\bibitem[Krummenacher et~al.(2016)Krummenacher, McWilliams, Kilcher, Buhmann,
  and Meinshausen]{conf/nips/KrummenacherMKB16}
Krummenacher, G., McWilliams, B., Kilcher, Y., Buhmann, J.~M., and Meinshausen,
  N.
\newblock Scalable adaptive stochastic optimization using random projections.
\newblock In \emph{NIPS}, pp.\  1750--1758, 2016.

\bibitem[Liberty et~al.(2007)Liberty, Woolfe, Martinsson, Rokhlin, and
  Tygert]{liberty2007randomized}
Liberty, E., Woolfe, F., Martinsson, P.G., Rokhlin, V., and Tygert, M.
\newblock Randomized algorithms for the low-rank approximation of matrices.
\newblock \emph{Proceedings of the National Academy of Sciences}, 104\penalty0
  (51):\penalty0 20167, 2007.

\bibitem[Lichman(2013)]{Lichman:2013}
Lichman, M.
\newblock {UCI} machine learning repository, 2013.
\newblock URL \url{http://archive.ics.uci.edu/ml}.

\bibitem[Luo et~al.(2016)Luo, Agarwal, Cesa-Bianchi, and
  Langford]{conf/nips/LuoACL16}
Luo, H., Agarwal, A., Cesa-Bianchi, N., and Langford, J.
\newblock Efficient second order online learning by sketching.
\newblock In \emph{NIPS}, pp.\  902--910, 2016.

\bibitem[Mahoney \& Drineas(2009)Mahoney and Drineas]{mahoney2009matrix}
Mahoney, M.~W. and Drineas, P.
\newblock Cur matrix decompositions for improved data analysis.
\newblock \emph{Proceedings of the National Academy of Sciences}, 106\penalty0
  (3):\penalty0 697--702, 2009.

\bibitem[Martinsson et~al.(2006)Martinsson, Rokhlin, and Tygert]{martinsson06}
Martinsson, P.-G., Rokhlin, V., and Tygert, M.
\newblock A randomized algorithm for the approximation of matrices.
\newblock \emph{Tech. Rep., Yale University, Department of Computer Science},
  \penalty0 (1361), 2006.

\bibitem[Melgaard \& Gu(2015)Melgaard and Gu]{journals/corr/MelgaardG15}
Melgaard, C. and Gu, M.
\newblock Gaussian elimination with randomized complete pivoting.
\newblock \emph{CoRR}, abs/1511.08528, 2015.

\bibitem[Miranian \& Gu(2003)Miranian and Gu]{miranian03}
Miranian, L. and Gu, M.
\newblock Stong rank revealing lu factorizations.
\newblock \emph{Linear Algebra and its Applications}, 367:\penalty0 1--16,
  2003.

\bibitem[NASA()]{nasa2jpg}
NASA.
\newblock Nasa celebrates 50 years of spacewalking.
\newblock URL
  \url{https://www.nasa.gov/image-feature/nasa-celebrates-50-years-of-spacewalking}.
\newblock Accessed on August 22, 2015. Published on June 3, 2015. Original
  photograph from February 7, 1984.

\bibitem[Pan(2000)]{pan00}
Pan, C.-T.
\newblock On the existence and computation of rank-revealing lu factorizations.
\newblock \emph{Linear Algebra and its Applications}, 316\penalty0
  (1):\penalty0 199--222, 2000.

\bibitem[Papadimitriou et~al.(2000)Papadimitriou, Raghavan, Tamaki, and
  Vempala]{journals/jcss/PapadimitriouRTV00}
Papadimitriou, C.~H., Raghavan, P., Tamaki, H., and Vempala, S.
\newblock Latent semantic indexing: A probabilistic analysis.
\newblock \emph{J. Comput. Syst. Sci.}, 61\penalty0 (2):\penalty0 217--235,
  2000.

\bibitem[Sarlos(2006)]{conf/focs/Sarlos06}
Sarlos, T.
\newblock Improved approximation algorithms for large matrices via random
  projections.
\newblock In \emph{FOCS}, pp.\  143--152. IEEE Computer Society, 2006.

\bibitem[Stewart(1992)]{journals/tsp/Stewart92}
Stewart, G.~W.
\newblock An updating algorithm for subspace tracking.
\newblock \emph{IEEE Trans. Signal Processing}, 40\penalty0 (6):\penalty0
  1535--1541, 1992.

\bibitem[Ubaru et~al.(2015)Ubaru, Mazumdar, and Saad]{conf/icml/UbaruMS15}
Ubaru, S., Mazumdar, A., and Saad, Y.
\newblock Low rank approximation using error correcting coding matrices.
\newblock In \emph{ICML}, volume~37, pp.\  702--710, 2015.

\bibitem[Wang et~al.(2016)Wang, Mehdad, Radev, and Stent]{conf/naacl/WangMRS16}
Wang, W.~Y., Mehdad, Y., Radev, D.~R., and Stent, A.
\newblock A low-rank approximation approach to learning joint embeddings of
  news stories and images for timeline summarization.
\newblock pp.\  58--68, 2016.

\bibitem[Woolfe et~al.(2008)Woolfe, Liberty, Rokhlin, and
  Tygert]{woolfe2008fast}
Woolfe, F., Liberty, E., Rokhlin, V., and Tygert, M.
\newblock A fast randomized algorithm for the approximation of matrices.
\newblock \emph{Applied and Computational Harmonic Analysis}, 25\penalty0
  (3):\penalty0 335--366, 2008.

\end{thebibliography}
\bibliographystyle{icml2017}

\onecolumn
\icmltitle{An Efficient, Sparsity-Preserving Online Algorithm for Data Approximation: \\ 
           Supplementary Material}



\vskip 0.3in


\section{The Singular Value Decomposition (SVD)}
For any real matrix $\mbf{A}\in\mathbb{R}^{m\times n}$ there exist orthogonal matrices $\mbf{U}\in\mathbb{R}^{m\times m}$ and $\mbf{V}\in\mathbb{R}^{n\times n}$ such that
\begin{eqnarray}
	\mbf{U}^T\mbf{A}\mbf{V} = \text{diag}\left(\sigma_1,\cdots,\sigma_p\right) \stackrel{\text{def}}{=} \Sigma\nonumber
	\end{eqnarray}
such that $p=\min(m,n)$ and $\sigma_1\ge \cdots \ge\sigma_p\ge 0$.  The decomposition $\mbf{A}=\mbf{U}\Sigma\mbf{V}^T$ is known as the Singular Value Decomposition~\cite{golub13}.

For a given matrix $\mathbf{A}$ with rank $\rho$ and a target rank $k$, rank-$k$ approximation using the SVD achieves the minimal residual error in both spectral and Frobenius norms:
\begin{customthm}{(Eckart-Young \cite{eckart1936approximation,golub13})}
\begin{align*}
\min_{rank\left(\mathbf{B}\right) \le k}\left|\left|\mathbf{A-B}\right|\right|_\xi^2 = 
\left|\left| \mathbf{A} - \mathbf{A}_k\right|\right|_\xi^2 = \sum_{j=k+1}^{\rho} \sigma_j\left(\mathbf{A}\right)
\end{align*}
where $\xi=F$ or $2$.
\end{customthm}

\section{Further Discussion of Rank-Revealing Algorithms}
An important class of algorithms against which we test SRLU is \textbf{rank-revealing} algorithms for low-rank approximation:
\begin{definition}
An LU factorization is \textbf{rank-revealing} \cite{miranian03} if
\begin{eqnarray}
\sigma_k\left(\mbf{A}\right)\ge\sigma_{\min}\left(\mbf{L}_{11}\mbf{U}_{11}\right)\gg\sigma_{\max}\left(\mbf{S}\right)\ge\sigma_{k+1}\left(\mbf{A}\right)\approx0.\nonumber
\end{eqnarray}
\end{definition}
Several drawbacks exists to the above definition, including that $\mbf{L}_{11}\mbf{U}_{11}$ is not a low-rank approximation of the original data matrix, and that only certain singular values are bounded.  Stronger algorithms were developed in \cite{miranian03} by modifying the definition above to create \text{strong rank-revealing} algorithms:
\begin{definition}
An LU factorization is \textbf{strong rank-revealing} if
\begin{enumerate}
\item
\begin{eqnarray}
\sigma_i\left(\mbf{A}_{11}\right) &\ge& \frac{\sigma_i\left(\mbf{A}\right)}{q_1\left(k,m,n\right)},\nonumber\\
\sigma_j\left(\mbf{S}\right) &\le& \sigma_{k+j}\left(\mbf{A}\right)q_1\left(k,m,n\right),\nonumber
\end{eqnarray}
\item
\begin{eqnarray}
\left|\left(\mbf{A}_{21}\mbf{A}_{11}^{-1}\right)|{ij}\right| &\le& q_2\left(k,n,m\right),\nonumber
\end{eqnarray}
\item
\begin{eqnarray}
\left|\left(\mbf{A}_{11}^{-1}\mbf{A}_{12}\right)|{ij}\right| &\le& q_3\left(k,n,m\right),\nonumber
\end{eqnarray}
\end{enumerate}
where $1\le i\le k,1\le j\le n-k$, and $q_1\left(k,m,n\right)$, $q_2\left(k,m,n\right),$ and $q_3\left(k,m,n\right)$ are functions bounded by low-degree polynomials of $k,m,$ and $n$.
\end{definition}
Strong rank-revealing algorithms bound all singular values of the submatrix $\mbf{A}_{11}$, but, as before, do not produce a low-rank approximation.  Furthermore, they require bounding approximations of the left and right null spaces of the data matrix, which is both costly and not strictly necessary for the creation of a low-rank approximation.  No known algorithms or numeric experiments demonstrate that strong rank-revealing algorithms can indeed be implemented efficiently in practice.

\section{Updating $\mbf{R}$}
The goal of TRLUCP is to access the entire matrix once in the initial random projection, and then choose column pivots at each iteration without accessing the Schur complement.  Therefore, a projection of the Schur complement must be obtained at each iteration without accessing the Schur complement, a method that first appeared in \cite{journals/corr/MelgaardG15}.  Assume that $s$ iterations of TRLUCP have been performed and denote the projection matrix
\[ \mbf{\Omega} = \bordermatrix{ & \text{\tiny $sb$} & \text{\tiny $b$} & \text{\tiny $n-(s+1)b$}  \cr & \mbf{\Omega}_1 & \mathbf{\Omega}_2 & \mbf{\Omega}_3}. \]
Then the current projection of the Schur complement is
\[ \mbf{R}^{\text{cur}} = \bordermatrix{ & \text{\tiny $b$} & \text{\tiny $n-(s+1)b$}  \cr & \mbf{R}_1^{\text{cur}} & \mathbf{R}_2^{\text{cur}}}=\begin{pmatrix}\mbf{\Omega}_2&\mbf{\Omega}_3\end{pmatrix}\begin{pmatrix}\mbf{S}_{11}&\mbf{S}_{12}\\\mbf{S}_{21}&\mbf{S}_{22}\end{pmatrix}, \]
where the right-most matrix is the current Schur complement.
The next iteration of TRLUCP will need to choose columns based on a random projection of the Schur complement, which we wish to avoid accessing.  We can write:

\begin{eqnarray}
\mbf{R}^{\text{update}}	&=& \mbf{\Omega}_3\left(\mbf{A}_{33}-\mbf{A}_{32}\mbf{A}_{22}^{-1}\mbf{A}_{23}\right)\nonumber\\
					&=& \mbf{\Omega}_3\mbf{A}_{33}+\mbf{\Omega}_2\mbf{A}_{23}-\mbf{\Omega}_2\mbf{A}_{23}-\mbf{\Omega}_3\mbf{A}_{32}\mbf{A}_{22}^{-1}\mbf{A}_{23}\nonumber\\
					&=& \mbf{\Omega}_3\mbf{A}_{33}+\mbf{\Omega}_2\mbf{A}_{23}-\mbf{\Omega}_2\mbf{L}_{22}\mbf{U}_{23}-\mbf{\Omega}_3\mbf{L}_{32}\mbf{U}_{23}\nonumber\\
					&=& \mbf{R}_2^{\text{current}}-\left(\mbf{\Omega}_2\mbf{L}_{22}+\mbf{\Omega}_3\mbf{L}_{32}\right)\mbf{U}_{23}.\label{update_fast}
\end{eqnarray}

Here the current $\mbf{L}$ and $\mbf{U}$ at stage $s$ have been blocked in the same way as $\Omega$.  Note equation (\ref{update_fast}) no longer has the term $\mbf{A}_{33}$.  Furthermore, $\mbf{A}_{22}^{-1}$ has been replaced by substituting in submatrices of $\mbf{L}$ and $\mbf{U}$ that have already been calculated, which helps eliminate potential instability.

When the block size $b = 1$ and TRLUCP runs fully ($k = {\bf min}(m,n)$), TRLUCP is mathematically equivalent to the Gaussian Elimination with Randomized Complete Pivoting (GERCP) algorithm of  \cite{journals/corr/MelgaardG15}. However, TRLUCP 
differs from GERCP in two very important aspects: TRLUCP is based on the Crout variant of the LU factorization, which allows efficient truncation for low-rank matrix approximation; and TRLUCP has been structured in block form for more efficient implementation.

\section{Proofs of Theorems}

\begin{customthm}{1}\label{theorem3}
For any truncated LU factorization
\begin{eqnarray}
\|\Pi_1\mathbf{A}\Pi_2^T-\lhat\uhat\| &=& \|\mathbf{S}\|\nonumber
\end{eqnarray}
for any norm $\|\cdot\|$.  Furthermore,
\begin{eqnarray}
\|\Pi_1\mathbf{A}\Pi_2^T-\left(\lhat\uhat\right)_s\|_2 &\le& 2\|\mathbf{S}\|_2+\sigma_{s+1}\left(\mathbf{A}\right)\nonumber
\end{eqnarray}
where $\left(\cdot\right)_s$ is the rank-$s$ truncated SVD for $s\le k\ll m,n$.
\end{customthm}
\begin{proof}
The equation simply follows from $\Pi_1\mathbf{A}\Pi_2^T = \lhat\uhat+\begin{pmatrix}0&0\\0&\mathbf{S}\end{pmatrix}$.  For the inequality:
\begin{eqnarray}
&&\|\Pi_1\mathbf{A}\Pi_2^T-\left(\lhat\uhat\right)_s\|_2 \nonumber\\
        &&\quad=\quad \|\Pi_1\mathbf{A}\Pi_2^T-\lhat\uhat+\lhat\uhat-\left(\lhat\uhat\right)_s\|_2\nonumber\\
	&&\quad\le\quad \|\Pi_1\mathbf{A}\Pi_2^T-\lhat\uhat\|_2+\|\lhat\uhat-\left(\lhat\uhat\right)_s\|_2\nonumber\\
	&&\quad=\quad \|\mathbf{S}\|_2+\sigma_{s+1}\left(\lhat\uhat\right)\nonumber\\
	&&\quad=\quad \|\mathbf{S}\|_2+\sigma_{s+1}\left(\Pi_1\mathbf{A}\Pi_2^T-\begin{pmatrix}0&0\\0&\mathbf{S}\end{pmatrix}\right)\nonumber\\
	&&\quad\le\quad \|\mathbf{S}\|_2 + \sigma_{s+1}\left(\mathbf{A}\right) + \|\mathbf{S}\|_2.\nonumber 
\end{eqnarray}
\end{proof}


\begin{customthm}{2}
For a general rank-$k$ truncated LU decomposition
\begin{eqnarray}
\sigma_j\left(\mathbf{A}\right) &\le& \sigma_j\left(\lhat\uhat\right)\left(1+\left(1+\frac{\|\mathbf{A}\|_2}{\sigma_k\left(\lhat\uhat\right)}\right)\frac{\|\mathbf{S}\|_2}{\sigma_j\left(\mathbf{A}\right)}\right).\nonumber
\end{eqnarray}
\end{customthm}
\begin{proof}
\begin{eqnarray}
&&\sigma_j\left(\mbf{A}\right)\nonumber\\
 	&&\quad\le\quad \sigma_j\left(\lhat\uhat\right)\left(1+\frac{\|\mbf{S}\|_2}{\sigma_j\left(\lhat\uhat\right)}\right)\nonumber\\
	&&\quad=\quad \sigma_j\left(\lhat\uhat\right)\left(1+\frac{\sigma_j\left(\mbf{A}\right)}{\sigma_j\left(\lhat\uhat\right)}\frac{\|\mbf{S}\|_2}{\sigma_j\left(\mbf{A}\right)}\right)\nonumber\\
	&&\quad\le\quad \sigma_j\left(\lhat\uhat\right)\left(1+\frac{\sigma_j\left(\lhat\uhat\right)+\|\mbf{S}\|_2}{\sigma_j\left(\lhat\uhat\right)}\frac{\|\mbf{S}\|_2}{\sigma_j\left(\mbf{A}\right)}\right)\nonumber\\
	&&\quad=\quad \sigma_j\left(\lhat\uhat\right)\left(1+\left(1+\frac{\|\mbf{S}\|_2}{\sigma_j\left(\lhat\uhat\right)}\right)\frac{\|\mbf{S}\|_2}{\sigma_j\left(\mbf{A}\right)}\right)\nonumber\\
	&&\quad\le\quad \sigma_j\left(\lhat\uhat\right)\left(1+\left(1+\frac{\|\mbf{S}\|_2}{\sigma_k\left(\lhat\uhat\right)}\right)\frac{\|\mbf{S}\|_2}{\sigma_j\left(\mbf{A}\right)}\right).\nonumber
\end{eqnarray}
\end{proof}
Note that the relaxation in the final step serves to establish a universal constant across all $j$, which leads to fewer terms that need bounding when the global SRLU swapping strategy is developed.


\begin{customthm}{3}\label{theorem8}
\begin{eqnarray}
\|\Pi_1\mbf{A}\Pi_2^T-\lhat\mbf{M}\uhat\|_2\phantom{.} &\le& 2\|\mbf{S}\|_2,\nonumber\\
\|\Pi_1\mbf{A}\Pi_2^T-\lhat\mbf{M}\uhat\|_F &\le& \|\mbf{S}\|_F.\nonumber
\end{eqnarray}
\end{customthm}
\begin{proof}
First
\begin{eqnarray}
&&\|\Pi_1\mbf{A}\Pi_2^T-\lhat\mbf{M}\uhat\|_2\nonumber\\
&&= \left|\left|\begin{pmatrix}0&\left(\mbf{Q}_1^L\right)^T\mbf{C}\left(\mbf{Q}_2^U\right)^T\\ \left(\mbf{Q}_2^L\right)^T\mbf{C}\left(\mbf{Q}_1^U\right)^T & \left(\mbf{Q}_2^L\right)^T\mbf{C}\left(\mbf{Q}_2^U\right)^T\end{pmatrix}\right|\right|_2\nonumber\\
&&\le \left|\left|\left(\mbf{Q}_1^L\right)^T\mbf{C}\left(\mbf{Q}_2^U\right)^T\right|\right|_2\nonumber\\
&&\quad+\left|\left|\begin{pmatrix}\left(\mbf{Q}_2^L\right)^T\mbf{C}\left(\mbf{Q}_1^U\right)^T & \left(\mbf{Q}_2^L\right)^T\mbf{C}\left(\mbf{Q}_2^U\right)^T\end{pmatrix}\right|\right|_2\nonumber\\
&&= \left|\left|\left(\mbf{Q}_1^L\right)^T\mbf{C}\left(\mbf{Q}_2^U\right)^T\right|\right|_2\nonumber\\
&&\quad+\left|\left|\left(\mbf{Q}_2^L\right)^T\mbf{C}\begin{pmatrix}\left(\mbf{Q}_1^U\right)^T&\left(\mbf{Q}_2^U\right)^T\end{pmatrix}\right|\right|_2\nonumber\\
&&\le 2\|\mbf{C}\|_2\nonumber\\
&&= 2\|\mbf{S}\|_2.\nonumber
\end{eqnarray}
Also
\begin{eqnarray}
&&\|\Pi_1\mbf{A}\Pi_2^T-\lhat\mbf{M}\uhat\|_F\nonumber\\
	&&= \left|\left|\begin{pmatrix}\mbf{Q}_1^L&\mbf{Q}_2^L\end{pmatrix}\begin{pmatrix}\left(\mbf{Q}_1^L\right)^T\\\left(\mbf{Q}_2^L\right)^T\end{pmatrix}\mbf{A}\begin{pmatrix}\left(\mbf{Q}_1^U\right)^T&\left(\mbf{Q}_2^U\right)^T\end{pmatrix}\begin{pmatrix}\mbf{Q}_1^U\\\mbf{Q}_2^U\end{pmatrix}-\mbf{Q}_1^L\left(\mbf{Q}_1^L\right)^T\mbf{A}\left(\mbf{Q}_1^U\right)^T\mbf{Q}_1^U\right|\right|_F\nonumber\\
	&&=\left|\left|\mbf{Q}_1^L\left(\mbf{Q}_1^L\right)^T\mbf{A}\left(\mbf{Q}_2^U\right)^T\mbf{Q}_2^U+\mbf{Q}_2^L\left(\mbf{Q}_2^L\right)^T\mbf{A}\left(\mbf{Q}_1^U\right)^T\mbf{Q}_1^U+\mbf{Q}_2^L\left(\mbf{Q}_2^L\right)^T\mbf{A}\left(\mbf{Q}_2^U\right)^T\mbf{Q}_2^U\right|\right|_F\nonumber\\
	&&=\left|\left|\mbf{Q}_1^L\left(\mbf{Q}_1^L\right)^T\mbf{C}\left(\mbf{Q}_2^U\right)^T\mbf{Q}_2^U+\mbf{Q}_2^L\left(\mbf{Q}_2^L\right)^T\mbf{C}\left(\mbf{Q}_1^U\right)^T\mbf{Q}_1^U+\mbf{Q}_2^L\left(\mbf{Q}_2^L\right)^T\mbf{C}\left(\mbf{Q}_2^U\right)^T\mbf{Q}_2^U\right|\right|_F\nonumber\\
	&&= \left|\left|\begin{pmatrix}\mbf{Q}_1^L&\mbf{Q}_2^L\end{pmatrix}\begin{pmatrix}0&\left(\mbf{Q}_1^L\right)^T\mbf{C}\left(\mbf{Q}_2^U\right)^T\\ \left(\mbf{Q}_2^L\right)^T\mbf{C}\left(\mbf{Q}_1^U\right)^T & \left(\mbf{Q}_2^L\right)^T\mbf{C}\left(\mbf{Q}_2^U\right)^T\end{pmatrix}\begin{pmatrix}\mbf{Q}_1^U\\\mbf{Q}_2^U\end{pmatrix}\right|\right|_F\nonumber\\
	&&= \left|\left|\begin{pmatrix}0&\left(\mbf{Q}_1^L\right)^T\mbf{C}\left(\mbf{Q}_2^U\right)^T\\ \left(\mbf{Q}_2^L\right)^T\mbf{C}\left(\mbf{Q}_1^U\right)^T & \left(\mbf{Q}_2^L\right)^T\mbf{C}\left(\mbf{Q}_2^U\right)^T\end{pmatrix}\right|\right|_F\nonumber\\
	&&\le \left|\left|\begin{pmatrix}\left(\mbf{Q}_1^L\right)^T\mbf{C}\left(\mbf{Q}_1^U\right)^T&\left(\mbf{Q}_1^L\right)^T\mbf{C}\left(\mbf{Q}_2^U\right)^T\\ \left(\mbf{Q}_2^L\right)^T\mbf{C}\left(\mbf{Q}_1^U\right)^T & \left(\mbf{Q}_2^L\right)^T\mbf{C}\left(\mbf{Q}_2^U\right)^T\end{pmatrix}\right|\right|_F\nonumber\\
	&&= \left|\left|\begin{pmatrix}\mbf{Q}_1^L&\mbf{Q}_2^L\end{pmatrix}\begin{pmatrix}\left(\mbf{Q}_1^L\right)^T\\\left(\mbf{Q}_2^L\right)^T\end{pmatrix}\mbf{C}\begin{pmatrix}\left(\mbf{Q}_1^U\right)^T&\left(\mbf{Q}_2^U\right)^T\end{pmatrix}\begin{pmatrix}\mbf{Q}_1^U\\\mbf{Q}_2^U\end{pmatrix}\right|\right|_F\nonumber\\
	&&= \|\mbf{C}\|_F\nonumber\\						
	&&= \|\mbf{S}\|_F.\nonumber
\end{eqnarray}
\end{proof}


\begin{customthm}{4}\label{theorem5}
SRP produces a rank-$k$ SRLU factorization with 
\begin{eqnarray}
\|\Pi_1\mbf{A}\Pi_2^T-\lhat\uhat\|_2 &\le& \gamma\sigma_{k+1}\left(\mbf{A}\right),\nonumber\\
\|\Pi_1\mbf{A}\Pi_2^T-\left(\lhat\uhat\right)_j\|_2 &\le& \sigma_{j+1}\left(\mbf{A}\right)\left(1+2\gamma\frac{\sigma_{k+1}\left(\mbf{A}\right)}{\sigma_j\left(\mbf{A}\right)}\right),\nonumber
\end{eqnarray}
where $j\le k$ and $\gamma=O\left(fk\sqrt{mn}\right)$.
\end{customthm}
\begin{proof}
Note that the definition of $\alpha$ implies
\begin{eqnarray}
\|\mbf{S}\|_2\le\sqrt{(m-k)(n-k)}|\alpha|.\nonumber
\end{eqnarray}
From \cite{pan00}:
\begin{eqnarray}
\sigma_{\min}\left(\overline{\mbf{A}}_{11}\right)\le\sigma_{k+1}\left(\mbf{A}\right).\nonumber
\end{eqnarray}
Then:
\begin{eqnarray}
\sigma_{k+1}^{-1}\left(\mbf{A}\right)	&\le& \|\overline{\mbf{A}}_{11}^{-1}\|_2\nonumber\\
							&\le& (k+1)\|\overline{\mbf{A}}_{11}^{-1}\|_{\max}\nonumber\\
							&\le& (k+1)\frac{f}{|\alpha|}.\nonumber
\end{eqnarray}
Thus
\begin{eqnarray}
|\alpha|\le f(k+1)\sigma_{k+1}\left(\mbf{A}\right).\nonumber
\end{eqnarray}
The theorem follows by using this result with Theorem \ref{theorem3}, with
\begin{eqnarray}
\gamma &\le& \sqrt{mn}f(k+1).\nonumber
\end{eqnarray}
\end{proof}


\begin{customthm}{5}
Assume the condition of SRLU (equation (2)) is satisfied.  Then for $1\le j\le k$:
\begin{eqnarray}
\frac{\sigma_j\left(\mbf{A}\right)}{1+\tau\frac{\sigma_{k+1}\left(\mbf{A}\right)}{\sigma_j\left(\mbf{A}\right)}} \le \sigma_j\left(\lhat\uhat\right) \le \sigma_j\left(\mbf{A}\right)\left(1+\tau\frac{\sigma_{k+1}\left(\mbf{A}\right)}{\sigma_j\left(\mbf{A}\right)}\right),\nonumber
\end{eqnarray}
where $\tau\le O\left(mnk^2f^3\right)$.
\end{customthm}
\begin{proof}
After running $k$ iterations of rank-revealing LU,
$$\Pi_1\mbf{A}\Pi_2^T=\lhat\uhat+\mbf{C},$$
where $\mbf{C}=\begin{pmatrix}0&0\\0&\mbf{S}\end{pmatrix}$, and $\mbf{S}$ is the Schur complement.  Then
\begin{eqnarray}
\sigma_j\left(\mbf{A}\right)	&\le& \sigma_j\left(\lhat\uhat\right)+\|\mbf{C}\|_2\nonumber\\
					&=& \sigma_j\left(\lhat\uhat\right)\left[1+\frac{\|\mbf{C}\|_2}{\sigma_j\left(\lhat\uhat\right)}\right].\label{eqn:detsing}
\end{eqnarray}
For the upper bound:
\begin{eqnarray}
\sigma_j\left(\lhat\uhat\right)	&=& \sigma_j\left(\mbf{A}-\mbf{C}\right)\nonumber\\
						&\le& \sigma_j\left(\mbf{A}\right)+\|\mbf{C}\|_2\nonumber\\
						&=& \sigma_j\left(\mbf{A}\right)\left[1+\frac{\|\mbf{C}\|_2}{\sigma_j\left(\mbf{A}\right)}\right]\nonumber\\
						&=& \sigma_j\left(\mbf{A}\right)\left[1+\frac{\|\mbf{S}\|_2}{\sigma_j\left(\mbf{A}\right)}\right].\nonumber
\end{eqnarray}
The final form is achieved using the same bound on $\gamma$ as in Theorem \ref{theorem5}.
\end{proof}


\begin{customthm}{6}
\begin{eqnarray}
\|\Pi_1\mbf{A}\Pi_2^T-\lhat\mbf{M}\uhat\|_2\phantom{.} &\le& 2\gamma\sigma_{k+1}\left(\mbf{A}\right),\nonumber\\
\|\Pi_1\mbf{A}\Pi_2^T-\lhat\mbf{M}\uhat\|_F &\le& \omega\sigma_{k+1}\left(\mbf{A}\right),\nonumber
\end{eqnarray}
where $\gamma=O\left(fk\sqrt{mn}\right)$ is the same as in Theorem \ref{theorem5}, and $\omega=O\left(fkmn\right)$.
\end{customthm}
\begin{proof}
Note that the definition of $\alpha$ implies
\begin{eqnarray}
\|\mbf{S}\|_F\le(m-k)(n-k)|\alpha|.\nonumber
\end{eqnarray}
The rest follows by using Theorem \ref{theorem8} in a manner similar to how Theorem \ref{theorem5} invoked Theorem \ref{theorem3}.
\end{proof}


\begin{customthm}{7}
If $\sigma_j^2\left(\mbf{A}\right)>2\|\mbf{S}\|_2^2$ then
\begin{eqnarray}
\sigma_j\left(\mbf{A}\right) \ge\sigma_j\left(\lhat\mbf{M}\uhat\right) \ge \sigma_j\left(\mbf{A}\right)\sqrt{1-2\gamma\left(\frac{\sigma_{k+1}\left(\mbf{A}\right)}{\sigma_j\left(\mbf{A}\right)}\right)^2},\nonumber
\end{eqnarray}
where $\gamma=O\left(mnk^2f^2\right)$, and $f$ is an input parameter controlling a tradeoff of quality vs. speed as before.
\end{customthm}
\begin{proof}
Perform QR and LQ decompositions $\lhat=\mbf{Q}_L\mbf{R}_L=:\begin{pmatrix}\mbf{Q}_1^L&\mbf{Q}_2^L\end{pmatrix}\begin{pmatrix}\mbf{R}_{11}^L&\mbf{R}_{12}^L\\ &\mbf{R}_{22}^L\end{pmatrix}$ and $\uhat=\mbf{L}_U\mbf{Q}_U=:\begin{pmatrix}\mbf{L}_{11}^U & \\ \mbf{L}_{21}^U&\mbf{L}_{22}^U\end{pmatrix}\begin{pmatrix}\mbf{Q}_1^U\\ \mbf{Q}_2^U\end{pmatrix}$.  Then
$$\lhat\mbf{M}\uhat=\mbf{Q}_1^L\left(\mbf{Q}_1^L\right)^T\mbf{A}\left(\mbf{Q}_1^U\right)^T\mbf{Q}_1^U.$$
Note that
\begin{eqnarray}
\mbf{A}^T\mbf{Q}_2^L	&=&\left(\lhat\uhat+\mbf{C}\right)^T\mbf{Q}_2^L\nonumber\\
					&=&\left(\mbf{Q}_1^L\mbf{R}_{11}^L\mbf{L}_{11}^U\mbf{Q}_1^U+\mbf{C}\right)^T\mbf{Q}_2^L\nonumber\\
					&=&\left(\mbf{Q}_1^U\right)^T\left(\mbf{L}_{11}^U\right)^T\left(\mbf{R}_{11}^L\right)^T\left(\mbf{Q}_1^L\right)^T\mbf{Q}_2^L+\mbf{C}^T\mbf{Q}_2^L\nonumber\\
					&=&\mbf{C}^T\mbf{Q}_2^L.\label{eqn:perpproj1}
\end{eqnarray}
Analogously
\begin{eqnarray}
\mbf{A}\left(\mbf{Q}_2^U\right)^T=\mbf{C}\left(\mbf{Q}_2^U\right)^T.\label{eqn:perpproj2}
\end{eqnarray}
Then
\footnotesize
\begin{eqnarray}
\sigma_j\left(\mbf{A}\right) &=& \sigma_j\begin{pmatrix}\left(\mbf{Q}_1^L\right)^T\mbf{A}\left(\mbf{Q}_1^U\right)^T&\left(\mbf{Q}_1^L\right)^T\mbf{A}\left(\mbf{Q}_2^U\right)^T\\ \left(\mbf{Q}_2^L\right)^T\mbf{A}\left(\mbf{Q}_1^U\right)^T & \left(\mbf{Q}_2^L\right)^T\mbf{A}\left(\mbf{Q}_2^U\right)^T\end{pmatrix}\nonumber\\
					&=& \sigma_j\begin{pmatrix}\left(\mbf{Q}_1^L\right)^T\mbf{A}\left(\mbf{Q}_1^U\right)^T&\left(\mbf{Q}_1^L\right)^T\mbf{C}\left(\mbf{Q}_2^U\right)^T\\ \left(\mbf{Q}_2^L\right)^T\mbf{C}\left(\mbf{Q}_1^U\right)^T & \left(\mbf{Q}_2^L\right)^T\mbf{C}\left(\mbf{Q}_2^U\right)^T\end{pmatrix}\nonumber\\
					&=&\sqrt{\lambda_j\left(\begin{pmatrix}\left(\mbf{Q}_1^L\right)^T\mbf{A}\left(\mbf{Q}_1^U\right)^T&\left(\mbf{Q}_1^L\right)^T\mbf{C}\left(\mbf{Q}_2^U\right)^T\\ \left(\mbf{Q}_2^L\right)^T\mbf{C}\left(\mbf{Q}_1^U\right)^T & \left(\mbf{Q}_2^L\right)^T\mbf{C}\left(\mbf{Q}_2^U\right)^T\end{pmatrix}^T\begin{pmatrix}\left(\mbf{Q}_1^L\right)^T\mbf{A}\left(\mbf{Q}_1^U\right)^T&\left(\mbf{Q}_1^L\right)^T\mbf{C}\left(\mbf{Q}_2^U\right)^T\\ \left(\mbf{Q}_2^L\right)^T\mbf{C}\left(\mbf{Q}_1^U\right)^T & \left(\mbf{Q}_2^L\right)^T\mbf{C}\left(\mbf{Q}_2^U\right)^T\end{pmatrix}\right)}\nonumber\\
					&=& \sqrt{\lambda_j\left(\begin{pmatrix}\left(\mbf{Q}_1^L\right)^T\mbf{A}\left(\mbf{Q}_1^U\right)^T&\left(\mbf{Q}_1^L\right)^T\mbf{C}\left(\mbf{Q}_2^U\right)^T\end{pmatrix}^T\begin{pmatrix}\left(\mbf{Q}_1^L\right)^T\mbf{A}\left(\mbf{Q}_1^U\right)^T&\left(\mbf{Q}_1^L\right)^T\mbf{C}\left(\mbf{Q}_2^U\right)^T\end{pmatrix}\color{white}\right)}\nonumber\\				
					&&\quad\quad\quad\quad\overline{\color{white}\left(\color{black}+\begin{pmatrix}\left(\mbf{Q}_2^L\right)^T\mbf{C}\left(\mbf{Q}_1^U\right)^T & \left(\mbf{Q}_2^L\right)^T\mbf{C}\left(\mbf{Q}_2^U\right)^T\end{pmatrix}^T\begin{pmatrix}\left(\mbf{Q}_2^L\right)^T\mbf{C}\left(\mbf{Q}_1^U\right)^T & \left(\mbf{Q}_2^L\right)^T\mbf{C}\left(\mbf{Q}_2^U\right)^T\end{pmatrix}\right)}\nonumber\\
					&\le& \sqrt{\lambda_j\left(\begin{pmatrix}\left(\mbf{Q}_1^L\right)^T\mbf{A}\left(\mbf{Q}_1^U\right)^T&\left(\mbf{Q}_1^L\right)^T\mbf{C}\left(\mbf{Q}_2^U\right)^T\end{pmatrix}^T\begin{pmatrix}\left(\mbf{Q}_1^L\right)^T\mbf{A}\left(\mbf{Q}_1^U\right)^T&\left(\mbf{Q}_1^L\right)^T\mbf{C}\left(\mbf{Q}_2^U\right)^T\end{pmatrix}\right)}\nonumber\\
					&&\quad\quad\quad\quad\overline{+\left|\left|\begin{pmatrix}\left(\mbf{Q}_2^L\right)^T\mbf{C}\left(\mbf{Q}_1^U\right)^T & \left(\mbf{Q}_2^L\right)^T\mbf{C}\left(\mbf{Q}_2^U\right)^T\end{pmatrix}\right|\right|_2}\nonumber\\
					&\le& \sqrt{\lambda_j\left(\begin{pmatrix}\left(\mbf{Q}_1^L\right)^T\mbf{A}\left(\mbf{Q}_1^U\right)^T&\left(\mbf{Q}_1^L\right)^T\mbf{C}\left(\mbf{Q}_2^U\right)^T\end{pmatrix}^T\begin{pmatrix}\left(\mbf{Q}_1^L\right)^T\mbf{A}\left(\mbf{Q}_1^U\right)^T&\left(\mbf{Q}_1^L\right)^T\mbf{C}\left(\mbf{Q}_2^U\right)^T\end{pmatrix}\right)+\|\mbf{C}\|_2^2}\nonumber\\
					&=& \sqrt{\lambda_j\left(\begin{pmatrix}\left(\mbf{Q}_1^L\right)^T\mbf{A}\left(\mbf{Q}_1^U\right)^T&\left(\mbf{Q}_1^L\right)^T\mbf{C}\left(\mbf{Q}_2^U\right)^T\end{pmatrix}\begin{pmatrix}\left(\mbf{Q}_1^L\right)^T\mbf{A}\left(\mbf{Q}_1^U\right)^T&\left(\mbf{Q}_1^L\right)^T\mbf{C}\left(\mbf{Q}_2^U\right)^T\end{pmatrix}^T\right)+\|\mbf{C}\|_2^2}\nonumber\\
					&=& \sqrt{\lambda_j\left(\left(\mbf{Q}_1^L\right)^T\mbf{A}\left(\mbf{Q}_1^U\right)^T\left(\left(\mbf{Q}_1^L\right)^T\mbf{A}\left(\mbf{Q}_1^U\right)^T\right)^T+\left(\mbf{Q}_1^L\right)^T\mbf{C}\left(\mbf{Q}_2^U\right)^T\left(\left(\mbf{Q}_1^L\right)^T\mbf{C}\left(\mbf{Q}_2^U\right)^T\right)^T\right)+\|\mbf{C}\|_2^2}\nonumber\\
					&\le& \sqrt{\lambda_j\left(\left(\mbf{Q}_1^L\right)^T\mbf{A}\left(\mbf{Q}_1^U\right)^T\left(\left(\mbf{Q}_1^L\right)^T\mbf{A}\left(\mbf{Q}_1^U\right)^T\right)^T\right)+\left|\left|\left(\mbf{Q}_1^L\right)^T\mbf{C}\left(\mbf{Q}_2^U\right)^T\left(\left(\mbf{Q}_1^L\right)^T\mbf{C}\left(\mbf{Q}_2^U\right)^T\right)^T\right|\right|_2+\|\mbf{C}\|_2^2}\nonumber\\
					&\le& \sqrt{\lambda_j\left(\left(\mbf{Q}_1^L\right)^T\mbf{A}\left(\mbf{Q}_1^U\right)^T\left(\left(\mbf{Q}_1^L\right)^T\mbf{A}\left(\mbf{Q}_1^U\right)^T\right)^T\right)+2\|\mbf{C}\|_2^2}\nonumber\\
					&\le& \sqrt{\sigma_j^2\left(\left(\mbf{Q}_1^L\right)^T\mbf{A}\left(\mbf{Q}_1^U\right)^T\right)+2\|\mbf{C}\|_2^2}\nonumber\\
					&=& \sqrt{\sigma_j^2\left(\lhat\mbf{M}\uhat\right)+2\|\mbf{C}\|_2^2}\nonumber\\
					&=& \sigma_j\left(\lhat\mbf{M}\uhat\right)\sqrt{1+2\left(\frac{\|\mbf{C}\|_2}{\sigma_j\left(\lhat\mbf{M}\uhat\right)}\right)^2}\nonumber\\
					&=& \sigma_j\left(\lhat\mbf{M}\uhat\right)\sqrt{1+2\left(\frac{\|\mbf{S}\|_2}{\sigma_j\left(\lhat\mbf{M}\uhat\right)}\right)^2}.\nonumber
\end{eqnarray}
\normalsize
Solve for $\sigma_j\left(\lhat\mbf{M}\uhat\right)$ for the lower bound.  The upper bound:
\begin{eqnarray}
\sigma_j\left(\mbf{A}\right) &=& \sigma_j\begin{pmatrix}\left(\mbf{Q}_1^L\right)^T\mbf{A}\left(\mbf{Q}_1^U\right)^T&\left(\mbf{Q}_1^L\right)^T\mbf{A}\left(\mbf{Q}_2^U\right)^T\\ \left(\mbf{Q}_2^L\right)^T\mbf{A}\left(\mbf{Q}_1^U\right)^T & \left(\mbf{Q}_2^L\right)^T\mbf{A}\left(\mbf{Q}_2^U\right)^T\end{pmatrix}\nonumber\\
	&\ge& \sigma_j\left(\left(\mbf{Q}_1^L\right)^T\mbf{A}\left(\mbf{Q}_1^U\right)^T\right)\nonumber\\
	&=& \sigma_j\left(\mathbf{Q}_1^L\left(\mbf{Q}_1^L\right)^T\mbf{A}\left(\mbf{Q}_1^U\right)^T\mathbf{Q}_1^U\right)\nonumber\\
	&=& \sigma_j\left(\lhat\mbf{M}\uhat\right).\nonumber
\end{eqnarray}
\end{proof}

\section{Analysis of the Choice of Block Size for SRLU}

A heuristic for choosing a block size for TRLUCP is described here, which differs from standard block size methodologies for the LU decomposition.  Note that a key difference of SRLU and TRLUCP from previous works is the size of the random projection: here the size is relative to the block size, not the target rank $k$ ($2pmn$ flops for TRLUCP versus the significantly larger $2kmn$ for others).  This also implies a change to the block size also changes the flop count, and, to our knowledge, this is the first algorithm where the choice of block size affects the flop count.  For problems where LAPACK chooses $b=64$, our experiments have shown block sizes of 8 to 20 to be optimal for TRLUCP.  Because the ideal block size depends on many parameters, such as the architecture of the computer and the costs for various arithmetic, logic, and memory operations, guidelines are sought instead of an exact determination of the most efficient block size.  To simplify calculations, only the matrix multiplication operations are considered, which are the bottleneck of computation.  Using standard communication-avoiding analysis, a good block size can be calculated with the following model: let $M$ denote the size of cache, $f$ and $m$ the number of flops and memory movements, and $t_f$ and $t_m$ the cost of a floating point operation and the cost of a memory movement.  We seek to choose a block size to minimize the total calculation time $T$ modeled as
\begin{eqnarray}
T &=& f\cdot t_f+m\cdot t_m.\nonumber
\end{eqnarray}
Choosing $p=b+c$ for a small, fixed constant $c$, and minimizing implies
\begin{eqnarray}
T &=& \left[(m+n-k)\left(k^2-kb\right)-\frac{4}{3}k^3+2bk^2-\frac{2}{3}b^2k\right]\cdot t_f\nonumber\\
   && + \left[(m+n-k)\left(\frac{k^2}{b}-k\right)-\frac{4}{3}\frac{k^3}{b}+2k^2-\frac{2}{3}bk\right] \cdot \frac{M}{\left(\sqrt{b^2+M}-b\right)^2}\cdot t_m.\nonumber
\end{eqnarray}
Given hardware-dependent parameters $M$, $t_f$, and $t_m$, a minimizing $b$ can easily be found.  

This result is derived as follows: we analyze blocking by allowing different block sizes in each dimension.  For matrices $\mbf{\Omega}\in\mathbb{R}^{p\times m}$ and $\mbf{A}\in\mathbb{R}^{m\times n}$ consider blocking in the form
\[
\mbf{\Omega}\cdot\mbf{R}=
\bordermatrix{ & \text{\tiny $\ell$} &  \text{\tiny $$}  & \text{\tiny $$} \cr \text{\tiny $s$}  & * & * & * \cr \text{\tiny $$}  & * & * & *}
\cdot
\bordermatrix{ & \text{\tiny $b$} &  \text{\tiny $$}  & \text{\tiny $$} \cr \text{\tiny $\ell$}  & * & * & * \cr \text{\tiny $$}  & * & * & * \cr \text{\tiny $$}  & * & * & *}.
\]
Then a current block update requires cache storage of
$$s\ell+\ell b+sb\le M.$$
Thus we will constrain
$$\ell\le\frac{M-sb}{s+b}.$$
The total runtime $T$ is
\begin{eqnarray}
T	&=& 2pmn\cdot t_f+\left(\frac{p}{s}\right)\left(\frac{m}{\ell}\right)\left(\frac{n}{b}\right)\left(s\ell+\ell b+sb\right)\cdot t_m\nonumber\\
	&=& 2pmn\cdot t_f+pmn\left(\frac{s+b}{sb}+\frac{1}{\ell}\right)\cdot t_m\nonumber\\
	&\ge& 2pmn\cdot t_f+pmn\left(\frac{s+b}{sb}+\frac{s+b}{M-sb}\right)\cdot t_m\nonumber\\
	&=& 2pmn\cdot t_f+pmnM\left(\frac{s+b}{sb\left(M-sb\right)}\right)\cdot t_m\nonumber\\
	&=:& 2pmn\cdot t_f+pmnML\left(s,b,M\right)\cdot t_m.\nonumber 
\end{eqnarray}
Given $\mbf{\Omega}$ and $\mbf{A}$, changing the block sizes has no effect on the flop count.  Optimizing $L\left(s,b,M\right)$ over $s$ yields
$$s^2+2sb=M.$$
By symmetry
$$b^2+2sb=M.$$
Note, nevertheless, that $s\le p$ by definition.  Hence
$$s^*=\min\left(\sqrt{\frac{M}{3}},p\right),$$
and
$$b^*=\max\left(\sqrt{\frac{M}{3}},\sqrt{p^2+M}-p\right).$$
These values assume
$$\ell^*=\frac{M-sb}{s+b}=\max\left(\sqrt{\frac{M}{3}},\sqrt{p^2+M}-p\right)=b^*.$$
This analysis applies to matrix-matrix multiplication where the matrices are fixed and the leading matrix is short and fat or the trailing matrix is tall and skinny.  As noted above, nevertheless, the oversampling parameter $p$ is a constant amount larger than the block size used during the LU factorization.  The total initialization time is
\begin{eqnarray}
T^{\text{init}}	&=& 2pmn\cdot t_f+pmnM\left(\frac{s+b}{sb\left(M-sb\right)}\right)\cdot t_m\nonumber\\
			&=& 2pmn\cdot t_f+mn\cdot\min\left(3\sqrt{3}\frac{p}{\sqrt{M}},\frac{M}{\left(\sqrt{p^2+M}-p\right)^2}\right)\cdot t_m.\nonumber
\end{eqnarray}
We next choose the parameter $b$ used for blocking the LU factorization, where $p=b+O\left(1\right)$.  The cumulative matrix multiplication ({\tt DGEMM}) runtime is

\begin{eqnarray}
T^{\text{\tt DGEMM}}	&=& \sum_{j=b:b:k-b}\left[2jb(m-j)+2jb(n-j-b)\right]\cdot t_f+2\left[j(m-j)+j(n-j-b)\right]\frac{M}{\left(\sqrt{b^2+M}-b\right)^2}\cdot t_m\nonumber\\
				&=& \left[(m+n-k)\left(k^2-kb\right)-\frac{4}{3}k^3+2bk^2-\frac{2}{3}b^2k\right]\cdot t_f+\nonumber\\
				&& + \left[(m+n-k)\left(\frac{k^2}{b}-k\right)-\frac{4}{3}\frac{k^3}{b}+2k^2-\frac{2}{3}bk\right]\frac{M}{\left(\sqrt{b^2+M}-b\right)^2}\cdot t_m\nonumber\\
				&=:& N_f^{\tt DGEMM}\cdot t_f+N_m^{\tt DGEMM}\cdot t_m.\nonumber
\end{eqnarray}

\begin{figure}[t]
\begin{center}
\centerline{\includegraphics[width=5in]{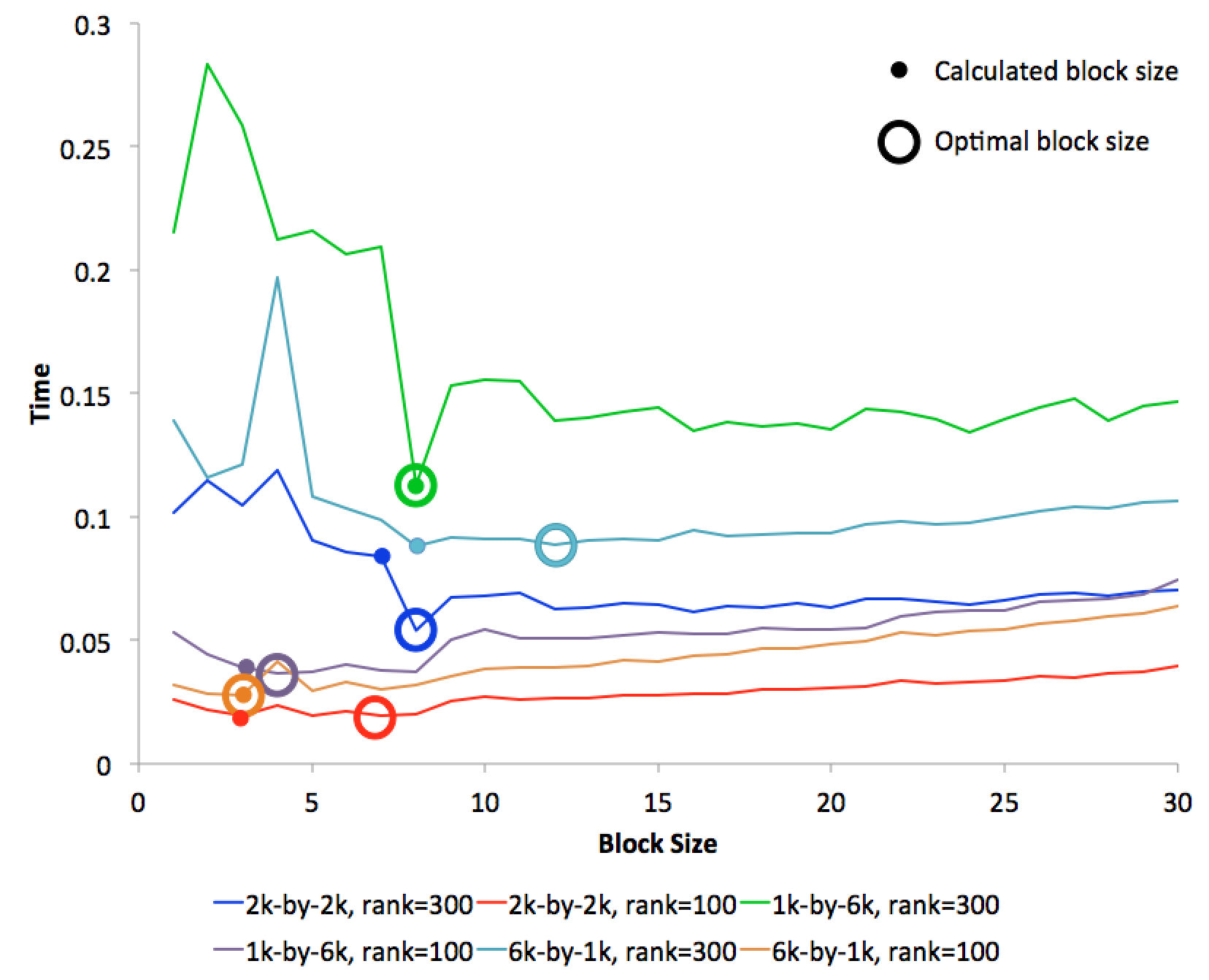}}
\caption{Benchmarking TRLUCP with various block sizes on random matrices of different sizes and truncation ranks.}
\label{block_size}
\end{center}
\vskip -0.2in
\end{figure}

The methodology for choosing a block size is compared to other choices of block size in Figure \ref{block_size}.  Note that LAPACK generally chooses a block size of 64 for these matrices, which is suboptimal in all cases, and can be up to twice as slow.  In all of the cases tested, the calculated block size is close to or exactly the optimal block size.

\section{Additional Notes and Experiments}

\subsection{Efficiency of SRLU}
Not only is the TRLUCP component efficient compared with other low-rank approximation algorithms, but also it becomes arbitrarily faster than the standard right-looking LU decomposition as the data size increases.  Because the LU decomposition is known to be efficient compared to algorithms such as the SVD \cite{demmel97}, comparing TRLUCP to right-looking LU exemplifies its efficiency, even though right-looking LU is not a low-rank approximation algorithm.

\begin{figure}[ht]
\begin{center}
\centerline{\includegraphics[width=3.8in]{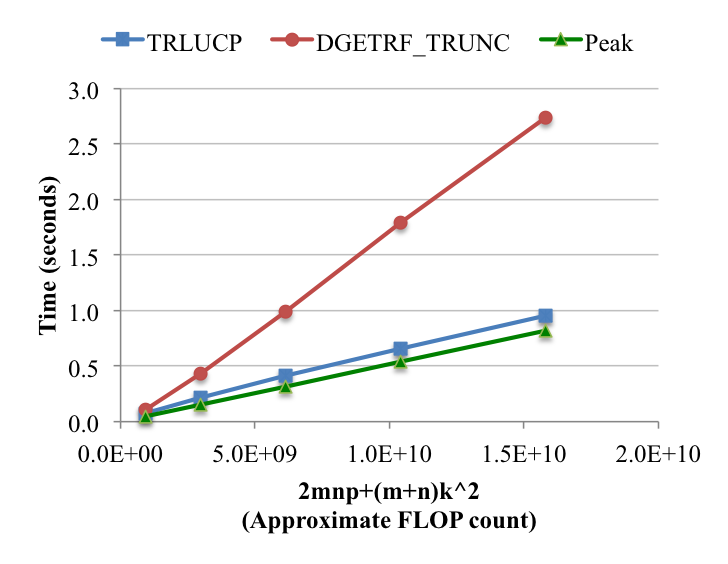}}
\caption{Computation time of TRLUCP versus the efficiency LU decomposition.}
\label{efficiency}
\end{center}
\vskip -0.2in
\end{figure} 

In Figure \ref{efficiency}, TRLUCP is benchmarked against truncated right-looking LU (called using a truncated version of the LAPACK library {\tt DGETRF}).  Experiments are run on random matrices, with the $x$-axis reflecting the approximate number of floating point operations.  Also plotted is the theoretical peak performance, which illustrates that TRLUCP is a highly efficient algorithm.




\subsection{Sparsity-Preservation}

Table \ref{sample-table} contains additional sparsity-preservation experiments on matrices from \cite{timdavis}.

\begin{table}[h]
  \caption{Sparsity preservation experiments of various sparse, non-symmetric data matrices.  The SRLU factorization is computed to 20\% of full-rank.  The Full SRLU factorization is the SRLU factorization with the Schur complement.  LU and SVD are the standard LU and SVD decompositions.  The SRLU relative error is the Frobenius-norm relative error of the SRLU factorization, which has a target rank that is 20 percent of the matrix rank.}
  \label{sample-table}
  \centering
  \begin{tabular}{llrlrrrrr}
    \toprule
    \multicolumn{2}{c}{Matrix Description}  & & & & \multicolumn{2}{c}{Nonzeros (rounded) In:} \\
    \cmidrule{1-3}\cmidrule{5-8}
    Name     			& Application     	& Nonzeros 	& & SRLU & Full SRLU & LU & SVD & SRLU Rel. Error\\
    \midrule
    \texttt{oscil\_dcop} 	& Circuits			& 1,544		& & 1,570	& 4.7K & 9.7K & 369K & 1.03e-3 \\
    \texttt{g7jac020} 		& Economics  		& 42,568 		& & 62.7K & 379K & 1.7M & 68M & 1.09e-6\\
    \texttt{tols1090} 		& Fluid dynamics	& 3,546		& & 2.2K & 4.7K & 4.6K & 2.2M & 1.18e-4\\
    \texttt{mhd1280a} 		& Electromagnetics	& 47,906		& & 184K & 831K & 129K & 3.3M & 4.98e-6\\
    \bottomrule
  \end{tabular}
\end{table}

\newpage
\subsection{Online Data Processing}

In many applications, reduced weight is given to old data.  In this context, multiplying the matrices $\mbf{U}_{11}$, $\mbf{U}_{12}$ and $\mbf{S}$ by some scaling factor less than 1 before applying spectrum-revealing pivoting will reflect the reduced importance of the old data.

\begin{figure}[h]
\centerline{\includegraphics[width=2.8in]{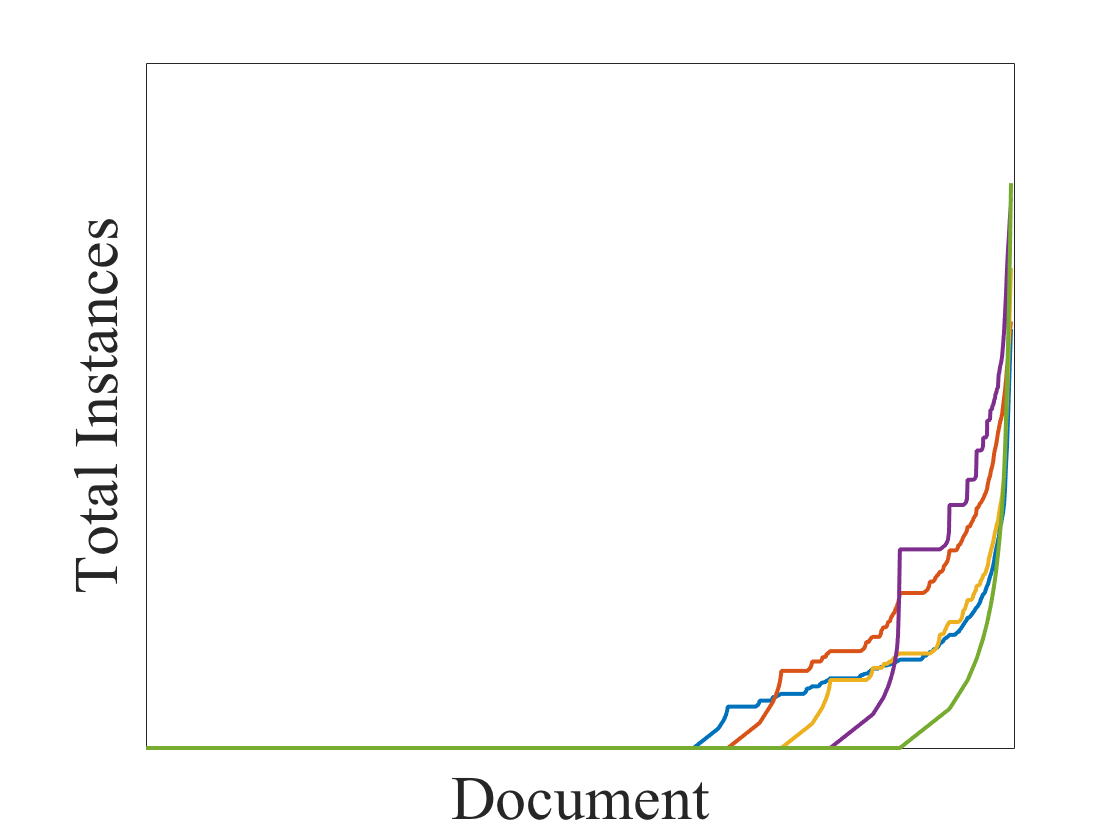}}
\caption{The cumulative uses of the top five most commonly used words in the Enron email corpus after reordering.}
\label{online}
\end{figure} 

The cumulative usages of the top 5 words in the Enron email corpus (after reordering) is plotted in Figure \ref{online}.  For the online updating experiment with the Enron email corpus, the covariance matrix of the top five most frequent words is

\centerline{\bordermatrix{ & \text{\tiny power} & \text{\tiny company} & \text{\tiny energy} & \text{\tiny market} & \text{\tiny california}
\cr \text{\tiny power} & 1 & 0.40 & 0.81 & 0.51 & 0.78
\cr \text{\tiny company} & 0.40 &1 & 0.42 & 0.57 & 0.28
\cr \text{\tiny energy} & 0.81 & 0.42 &1 & 0.51 & 0.78
\cr \text{\tiny market} & 0.51 & 0.57 & 0.51 & 1& 0.48
\cr \text{\tiny california} & 0.78 & 0.23 & 0.78 & 0.48 &1}.}




\end{document}